\documentclass[letterpaper,onecolumn,12pt]{article}  % Enable this line and disable the 
\usepackage{amssymb}
\usepackage{amsmath}
\usepackage{graphicx}  
\usepackage{epstopdf}
\usepackage{hyperref}
\usepackage{color}

\newtheorem{theorem}{Theorem}
\newtheorem{corollary}{Corollary}
\newtheorem{lemma}{Lemma}
\newtheorem{remark}{Remark}
\newtheorem{assumption}{Assumption}
\newtheorem{proposition}{Proposition}
\newtheorem{definition}{Definition}
\newtheorem{example}{Example}

\newcommand{\defeq}{\overset{\mathrm{def}}{=}}
\newcommand{\tr}{\intercal}
\newcommand{\pf}{\textit{Proof. }}
\newcommand{\qed}{\hfill$\diamond$}

\title{\textbf{Ambiguity Tube MPC}}

\author{Fan Wu, Mario E.~Villanueva, Boris Houska}
\date{ShanghaiTech University}

\begin{document}

\maketitle

% \begin{keyword}
% stochastic control, model predictive control, stability analysis, Markov processes, martingale theory
% \end{keyword}

\begin{abstract}
This paper is about a class of distributionally robust model predictive controllers (MPC) for nonlinear stochastic processes that evaluate risk and control performance measures by propagating ambiguity sets in the space of state probability measures. A framework for formulating such ambiguity tube MPC controllers is presented, which is based on modern measure-theoretic methods from the field of optimal transport theory. Moreover, a supermartingale based  analysis technique is proposed, leading to stochastic stability results for a large class of distributionally robust controllers for linear and nonlinear systems. In this context, we also discuss how to construct terminal cost functions for stochastic and distributionally robust MPC that ensure closed-loop stability and asymptotic convergence to robust invariant sets. The corresponding theoretical developments are illustrated by tutorial-style examples and a numerical case study.
\end{abstract}

\section{Introduction}
\label{sec::introduction}
Traditional robust MPC formulations that systematically take model uncertainties and external disturbances into account can be categorized into two classes. The first class of robust MPC controllers are based on min-max~\cite{Houska2019,Rawlings2009} or tube-based MPC formulations~\cite{Langson2004,Mayne2005,Rakovic2012}, which typically assume that worst-case bounds on the uncertainty are available. This is in contrast to the second class of optimization based robust controllers, namely stochastic MPC controllers~\cite{Kouvaritakis2015,Mesbah2016}, which assume that the probability distribution of the external disturbance is known. The main practical difference between these formulations is that most stochastic MPC controllers attempt to either bound or penalize the probability of a constraint violation, but, in contrast to min-max MPC formulations, conservative worst-case constraints are not enforced.

In terms of recent developments in the field of robust MPC, several attempts have been made to unify the above classes by considering distributionally robust MPC controllers~\cite{Parys2016}. Here, one assumes that the uncertainty is stochastic, but the associated probability distribution is only known to be in a given ambiguity set. Thus, in the most general setting, distributionally robust MPC formulations contain both traditional stochastic MPC as well as min-max MPC as special cases: in the context of stochastic MPC the ambiguity set is a singleton whereas min-max MPC is based on ambiguity sets that contain all uncertainty distributions with a given bounded support. Notice that modern distributionally robust MPC formulations are often formulated by using risk measures~\cite{Sopasakis2019}. This trend is motivated by the availability of rather general classes of \mbox{coherent---and}, most importantly, computationally \mbox{tractable---risk} measures, such as the conditional value at risk~\cite{Rockafellar2013}.

The current paper focuses on distributionally robust MPC problems that are formulated as ambiguity controllers. Here, the main idea is to propagate sets in the space of probability measures on the state space. The primary motivation for analyzing such classes of controllers is, however, not to develop yet another way to formulate robust MPC, but to develop a coherent stability theory for a very general class of distributionally robust MPC controllers that contain traditional tube based MPC as well as stochastic MPC as a special case.

Now, before we can outline why such a general stability theory for ambiguity controllers is of fundamental interest for control theory \mbox{research---especially}, in the emerging era of learning based MPC~\cite{Hewing2020,Zanon2021}, where uncertain models are omni-\mbox{present---one} has to first mention that there exist already many stability results for MPC. First of all, the stability of classical (certainty-equivalent) MPC has been analyzed by many \mbox{authors---be} it for tracking or economic MPC, with or without terminal costs or regions~\cite{Chen1998,Gruene2009,Rawlings2009}. Similarly, the stability of variants of min-max MPC schemes have been analyzed exhaustively~\cite{Mayne2005,Villanueva2017}, although the development of a unified stability analysis for general set-based MPC controllers is a topic of ongoing research~\cite{Villanueva2020b}.

The above reviewed results on the stability of classical certainty-equivalent and tube MPC controllers have in common that they are all based on the construction of Lyapunov functions, which descend along the closed-loop trajectories of the robust controller. This is in contrast to existing results on the stability of stochastic MPC, which are usually based on the theory of non-negative supermartingales~\cite{Doob1953,Feller1971}. The corresponding mathematical foundation for such stability results has been developed by H.J.~Kushner and R.S.~Bucy for analyzing general Markov processes. The corresponding original work can be found in~\cite{Bucy1965} and~\cite{Kushner1965}. We additionally refer to~\cite{Kushner2014} for an review of the history of this theory. At the current status of research on stochastic MPC, martingale theory has been applied to special classes of linear MPC controllers with multiplicative uncertainty~\cite{Bernardini2012}. Moreover, an impressive collection of articles by M.~Cannon and B.~Kouvaritakis has appeared during the last two decades, which has had significant impact on shaping the state-of-the-art of stochastic MPC. As we cannot possibly list all of the papers of these authors, we refer at this point to their textbook~\cite{Kouvaritakis2015} for an overview of formulations and stability results for stochastic linear MPC. Additionally, the book chapter~\cite{Kouvaritakis2016} comes along with an excellent overview of recursive feasibility results for stochastic MPC as well as a proof of stochastic stability with respect to ellipsoidal regions that are derived by using non-negative supermartingales, too. 

Given the above list of articles it can certainly be stated that the question how to establish stability results for both robust mix-max as well as stochastic MPC has received significant attention and much progress has been made. Nevertheless, looking back at the MPC literature from the last decade, it must also be stated that this question has raised a critical discussion. For example, a general critique of robust MPC can be found in~\cite{Mayne2015}, which points out the lack of a satisfying treatment of stabilizing terminal conditions for stochastic MPC. Similarly,~\cite{Chatterjee2015} points out various discrepancies between deterministic and stochastic MPC. From these articles, it certainly becomes clear that one has to distinguish carefully between rigorous supermartingale based stability analysis in the above reviewed sense of Kushner and Bucy and weaker properties of stochastic MPC from which stability may not be inferred. Among these weaker properties are bounds on the asymptotic average performance of stochastic MPC~\cite{Kouvaritakis2015,Kouvaritakis2016}, which do not necessarily imply stability. Moreover, in the past 5 years several articles have appeared, which all exploit input-to-state-stability assumptions for establishing convergence of stochastic MPC controllers to robust invariant sets~\cite{Lorenzen2016,Sehr2018}. One of the strongest results in this context appeared only a few weeks ago in~\cite{Munoz2020}, where an input-to-state-stability assumption in combination with the Borel-Cantelli lemma is used to establish conditions under which the state of a potentially nonlinear Markov process converges almost surely to a minimal robust invariant set. These conditions are applicable for establishing convergence of a variety of stochastic MPC formulations. Nevertheless, it has to be recalled here that, in general, neither stability implies convergence nor convergence implies stability. As such, none of the above reviewed articles proposes a completely satisfying answer to the question how asymptotic stability conditions can be established for general stochastic, let alone distributionally robust, MPC.

\paragraph*{Contribution.} This paper is about the mathematical formulation and stochastic stability analysis of distributionally robust MPC controllers for general, potentially nonlinear, but Lipschitz continuous stochastic discrete-time systems. Here, the focus is on so-called ambiguity tube MPC controllers that are based on the propagation of sets in the space of state probability measures, such that all theoretical results are applicable to traditional stochastic as well as tube based min-max MPC as special cases. The corresponding contributions of the current article can be outlined as follows.

\begin{enumerate}

\item Section~\ref{sec::AmbiguityControl} develops a novel framework for formulating ambiguity tube MPC problems by exploiting measure-theoretic concepts from the field of modern optimal transport theory~\cite{Villani2005}. In detail, we propose a Wasserstein metric based setting, which leads to a well-formulated class of ambiguity tube MPC controllers that admit a continuous value function; see Theorem~\ref{thm::existence}.\\

\item Section~\ref{sec::stability} presents a complete stability analysis for ambiguity tube MPC that is applicable to general Lipschitz continuous stochastic discrete-time system under mild assumptions on the coherentness of the optimized performance and risk measures, as well as on the consistency of the terminal cost function of the MPC controller. In detail, Theorem~\ref{thm::martingale} establishes conditions under which the cost function of the ambiguity tube MPC controller is a non-negative supermartingale along its closed-loop trajectories, which can then be used to establish robust stability \mbox{or---under} a slightly stronger regularity \mbox{assumption---robust} asymptotic stability of the closed loop system in a stochastic sense, as summarized in Theorems~\ref{thm::stability} and~\ref{thm::asympstability}. Notice that these stability results are both more general than the existing stability and convergence statements about stochastic MPC in~\cite{Chatterjee2015,Kouvaritakis2016,Munoz2020}, as they apply to general nonlinear systems and ambiguity set based formulations. Besides, Theorem~\ref{thm::asympstability} establishes conditions under which the closed-loop trajectories of stochastic (or ambiguity based) MPC is asymptotically stable with respect to a minimal robust invariant set. This is in contrast to the less tight results in~\cite{Kouvaritakis2016}, which only establish stability and convergence of linear stochastic MPC closed-loop systems with respect to an ellipsoidal set that overestimates the actual (in general, non-ellipsoidal) limit set of the stochastic ancillary closed-loop system.\\

\item Section~\ref{sec::implementation} discusses how to implement ambiguity tube MPC in practice. Here, our focus \mbox{is---for} the sake of simplicity of \mbox{presentation---on} linear systems, although remarks on how this can be implemented for nonlinear systems are provided, too. The purpose of this section is to illustrate how the technical assumptions from Section~\ref{sec::stability} can be satisfied in practice. In this context, a relevant technical contribution is summarized in Lemma~\ref{lem::assumptions}, which explains how to construct stabilizing terminal cost functions for stochastic and ambiguity tube MPC.
\end{enumerate}

Notice that as much this paper attempts to make a step forward towards a more coherent stability analysis and treatment of stabilizing terminal conditions for stochastic and distributionally robust MPC, it must be stated clearly that we do not claim to be anywhere close to addressing the long list of conceptual issues of robust MPC that D.~Mayne summarized in his critique~\cite{Mayne2015}. Nevertheless, in order to assess the role of this paper in the ongoing development of robust MPC, Section~\ref{sec::conclusions} does not only summarize and interpret the contributions of this paper, but also comments on the long list of open issues that research on robust MPC is currently facing.

\paragraph*{Notation.}
If $(R,r)$ is a metric space with given distance function $r: R \times R \to \mathbb R_+$, we use the notation $\mathbb K(R)$ to denote the set of compact subsets of $R$---assuming that it is clear from the context what $r$ is. Similarly, if $(R_1,r_1)$ and $(R_2,r_2)$ are two metric spaces, we use the notation $\mathcal L(R_1,R_2)$ to denote the set of Lipschitz continuous functions from~$R_1$ to~$R_2$ with respect to the distance functions $r_1$ and $r_2$. Moreover, $\mathcal L_1(R_1,R_2)$ denotes the subset of $\mathcal L(R_1,R_2)$ that consist of all functions from $R_1$ to $R_2$ whose Lipschitz constant is smaller than or equal to~$1$. Finally, $R_1 \times R_2$ is again a metric space with distance function $$r((a,b),(c,d)) \; \defeq \; r_1(a,c)+r_2(b,d)$$
for all $a,c \in R_1$ and all $b,d \in R_2$. If nothing else is stated, we assume that the distance function in the new metric space $R_1 \times R_2$ is constructed in the above additive \mbox{way---without} always saying this explicitly. Additionally, at some points in this paper, we make use of the notation
\[
\mathrm{dist}_q(a,B) \; \defeq \; \min_{b \in B} \; \| a-b \|_q
\]
to denote the distance of a point $a \in \mathbb R^n$ to a compact set $B \in \mathbb K(\mathbb R^n)$, where $\Vert \cdot \Vert_q$ denotes the H\"older $q$-norm.

For the special case that $R_1 \in \mathbb K(\mathbb R^n)$ and $R_2 \in \mathbb K(\mathbb R^m)$ are equipped with the standard Euclidean distance, we denote with
$$d_{\mathcal L}: \mathcal L(R_1,R_2) \times \mathcal L(R_1,R_2) \to \mathbb R_+$$
the corresponding Hilbert-Sobolev distance,
\[
d_{\mathcal L}(\mu,\nu) \; \defeq \; \sqrt{\int \| \mu-\nu \|_2^2 + \| \nabla \mu - \nabla  \nu \|_2^2 \; \mathrm{d}x} \; ,
\]
which is defined for all $\mu,\nu \in \mathcal L(R_1,R_2)$, where $\nabla$ denotes the weak gradient operator recalling that Lipschitz continuous functions are differentiable almost everywhere. With this notation, $(\mathcal L(R_1,R_2), d_{\mathcal L})$ is a metric space.

Next, let $X \in \mathbb K(\mathbb R^n)$ denote a given compact set in $\mathbb R^n$. We use the notation $\mathcal P(X)$ to denote the set of Borel probability measures on $X$. Notice that this definition is such that $p(X) = 1$ for all $p \in \mathcal P(X)$; and $\mathcal P(X)$ is convex. In this context, we additionally introduce the notation $\mathcal B(X)$ to denote the associated Borel $\sigma$-algebra of $X$. Notice that $\mathcal P(X)$ turns out to be a metric space with respect to the Wasserstein distance function\footnote{A very detailed review of the history and mathematical properties of Wasserstein distances can be found in~\cite[Chapter~6]{Villani2005}}, which is defined as follows.

\begin{definition}
Throughout this paper, we use the notation $d_{\mathrm{W}}: \mathcal P(X) \times \mathcal P(X) \to \mathbb R_+$ to denote the Wasserstein distance, given by
\[
\forall p,q \in \mathcal P(X), \quad d_{\mathrm{W}}(p,q) \; \defeq \sup_{\varphi \in \mathcal L_1(X,\mathbb R)} \left| \int \varphi \, \mathrm{d}(p-q) \right| \, .
\]
\end{definition}

Notice that $d_{\mathrm{W}}$ is well-defined in our context, as all Lipschitz continuous functions are $\mathcal B(X)$-measurable and, consequently, the integrals in the above definition exist and are finite, since we assume that $X$ is compact.

Throughout this paper, we use the notation $\delta_y \in \mathcal P(X)$ to denote the Dirac measure at a point $y \in X$, given by
\[
\forall Y \in \mathcal B(X), \qquad \delta_y(Y) \; \defeq \; \left\{
\begin{array}{ll}
1 \quad \text{if} \; y \in Y \\
0 \quad \text{otherwise} \; .
\end{array}
\right.
\]
Because this paper makes intense use of the concept of ambiguity sets, we additionally introduce the shorthand notation
$$\mathcal K(X) \; \defeq \; \mathbb K( \mathcal P(X))$$
to denote the set of subsets of $\mathcal P(X)$ that are compact in the Wasserstein space $(\mathcal P(X),d_{\mathrm{W}})$. By construction, $\mathcal K(X)$ is itself a metric space with respect to the Hausdorff-Wasserstein distance,
\[
\hspace{-0.2cm}
d_\mathrm{H}(P,Q) \defeq \max \left\{ \max_{p \in P} \min_{q \in Q} \, d_{\mathrm{W}}(p,q), \, \max_{q \in Q} \min_{p \in P} \, d_{\mathrm{W}} (p,q) \right\},
\]
defined for all $P,Q \in \mathcal K(X)$.

\begin{remark}
If $w$ is a random variable with given Lebesgue integrable probability distribution $\rho: \mathbb R^{n} \to \mathbb R_+$, the probability of the event $w \in W$ for a Borel set $W \subseteq \mathbb R^n$ is denoted by
\[
\mathrm{Pr}(w \in W) \ \defeq \ \int_W 1 \, \mathrm{d}p \ \defeq \ p(W) \ \defeq \ \int_{W} \rho(w) \, \mathrm{d}w,
\]
where $p \in \mathcal P(\mathbb R^{n})$ is called the probability measure of $w$. Notice that all four notations are, by definition, equivalent. Many articles on stochastic control, for instance,~\cite{Bucy1965,Kushner1965,Kushner2014} work with measures rather than probability distributions, as this has many technical advantages~\cite{Villani2005}. This means that we specify the probability measure $p$ rather than the probability distribution $\rho$. Notice that the relation
\[
\rho = \frac{\mathrm{d}p}{\mathrm{d}w}
\]
holds, where the right-hand expression denotes the Radon-Nikodyn derivative of the measure $p$ with respect to the traditional Lebesgue measure~\cite{Taylor1996}.
\end{remark}

%%%%%%%%%%%%%%%%%%%%%%%%%%%%%%%%%%%%%%%%%%%%%%%%%%%%

\section{Control of Ambiguity Tubes}
\label{sec::AmbiguityControl}

This section introduces a general class of uncertain control system models and their associated Markovian kernels that can be used to propagate state probability measures. Section~\ref{sec::AmbiguityTubes} exploits the properties of these Markovian kernels to introduce a topologically coherent framework for defining associated ambiguity tubes that have compact cross-section in $\mathcal K(X)$. Moreover, Section~\ref{sec::ambigiutyMeasures} focuses on an axiomatic characterization of proper risk and performance measures, which are then used in Section~\ref{sec::ambiguityMPC} to introduce a general class of ambiguity tube MPC controllers, completing our problem formulation.

\subsection{Uncertain Control Systems}
\label{sec::system}
This paper concerns uncertain nonlinear discrete-time control systems of the form
\begin{eqnarray}
\label{eq::system}
\forall k \in \mathbb N_0, \qquad x_{k+1} = f(x_k,u_k,w_k) \; ,
\end{eqnarray}
where $x_k$ denotes the state at time $k \in \mathbb N_0$, evolving in the state-space domain \mbox{$X \in \mathbb K(\mathbb R^{n_x})$}. Moreover, \mbox{$u_k \in U$} denotes the control input at time $k$ with domain \mbox{$U \in \mathbb K(\mathbb R^{n_u})$} and, finally, $w_k$ denotes an uncertain external disturbance with compact support \mbox{$W \in \mathbb K(\mathbb R^{n_w})$}.
\begin{assumption}
\label{ass::fLip}
We assume that the right-hand side function $f$ satisfies
\begin{eqnarray}
\label{eq::fLip}
f \in \mathcal L \left( \, X \times U \times W , \, X \, \right) \; .
\end{eqnarray}
\end{assumption}
The above assumption does not only require the function $f$ to be Lipschitz continuous, but it also requires that the image set of $f$ must be contained in $X$. Consequently, it should be mentioned that, throughout this paper, we interpret the set $X$ as a sufficiently large region of interest in which we wish to analyze the system's behavior. Consequently, if $f$ is Lipschitz continuous on $\mathbb R^{n_x}$ but not bounded, we redefine $f \leftarrow \mathrm{proj}_X \circ f$ with $\mathrm{proj}_X$ denoting a Lipschitz continuous projection onto $X$, such that~\eqref{eq::fLip} holds by construction. Notice that the set $X$ should not be mixed up with the set $\mathbb X \subseteq X$ that could, for example, model a so-called state-constraint; that is, a region in which we would like to keep the system state with high probability.

In the following, we additionally introduce a compact set $\mathcal U \in \mathbb K( \mathcal L(X,U) )$ in order to denote a class of ancillary feedback laws in the metric space $( \mathcal L(X,U), d_{\mathcal L} )$. In the context of this paper, $\mathcal U$ models a suitable class of computer representable feedback laws. 

\begin{example}
One of the easiest examples for a class of computer representable feedback policies is the set
\[
\mathcal U = \left\{ \, x \to Kx + k \, \middle| 
\begin{array}{l}
\exists k \in \mathbb R^{n_u}, \exists K \in \mathbb R^{n_u \times n_x}: \\
\ \ \| K \| \leq L 
\end{array}
 \right\}
\]
of affine control laws with bounded feedback gain, where $L < \infty$ is a given bound on the norm of the matrix $K$, such that all functions in $\mathcal U$ are Lipschitz continuous. Specific feedback laws can in this case be represented by storing the finite dimensional matrix $K$ and the offset vector $k$.
\end{example}

\subsection{Models of Stochastic Uncertainties}
\label{sec::stochastic}

In order to refine the above uncertain system model, we introduce the probability spaces $(W,\mathcal B(W),\omega_k)$. Here, $\omega_k \in \mathcal P(W)$ denotes the probability measures of the random variable \mbox{$w_k: W \to \mathbb R$} such that
\[
\forall W' \in \mathcal B(W), \quad \mathrm{Pr}( w_k \in W' ) = \omega_k(W') \; .
\]
In the most general setting, we might not know the measure $\omega_k$ up to a high precision, but we work with the more realistic assumption that an ambiguity set $\Omega \in \mathcal K(W)$ is given, which means that the sequence $\omega_k$ is only known to satisfy $\omega_k \in \Omega$ for all $k \in \mathbb N_0$.

In order to proceed with this modeling assumption, we analyze the closed-loop system
\begin{eqnarray}
\label{eq::ClosedLoopSystem}
\forall k \in \mathbb N_0, \qquad x_{k+1} = f(x_k,\mu(x_k),w_k) \; ,
\end{eqnarray}
for a given feedback law $\mu \in \mathcal U$. Because the disturbance sequence $w_0,w_1,\ldots $ consists of independent random variables, the states $x_k$ are random variables, too. Now, if \mbox{$p_k \in \mathcal P(X)$} denotes a probability measure that is associated with $x_k$, then the probability measure \mbox{$p_{k+1} \in \mathcal P(X)$} of the random variable $x_{k+1}$ in~\eqref{eq::ClosedLoopSystem} is a function of $p_k$, $\mu$, and $\omega_k$. Formally, this propagation of measures can be defined by using a parametric Markovian kernel, $\mathcal N[x,\mu,\omega]: \mathcal B(X) \to \mathbb R$,  given by
\[
\mathcal N[x,\mu,\omega](X^+) \defeq \omega \left( \; \left\{ \, w \in W \, \middle| \, f(x,\mu(x),w) \in X^+ \, \right\} \; \right)
\]
for all Borel sets $X^+ \in \mathcal B(X)$. The corresponding transition map
\begin{eqnarray}
\label{eq::PHI}
\Phi( p, \mu, \omega ) \; \defeq \; \int_X \mathcal N[x,\mu,\omega] \, p(\mathrm{d}x)
\end{eqnarray}
is then well-defined for all $p \in \mathcal P(X)$, all $\mu \in \mathcal U$, and all probability measures $\omega \in \mathcal P(W)$, where $\mathrm{d}x$ denotes the traditional Lebesgue measure~\cite{Feller1971}. This follows from our assumption that $f$ is Lipschitz continuous such that the Markovian kernel $\mathcal N[\cdot, \mu,\omega](X^+)$ is for any given Borel set $X^+$ a $\mathcal B(X)$-measurable function in $x$. In summary, the associated recursion for the sequence of measures $p_0,p_1,p_2, \ldots$ can be written in the form
\begin{eqnarray}
\label{eq::ProbSystem}
\forall k \in \mathbb N_0, \qquad p_{k+1} = \Phi( p_k, \mu, \omega_{k} ) \; .
\end{eqnarray}
The following technical lemma establishes an important property of the function $\Phi$.

\begin{lemma}
\label{lem::Lipschitz}
If Assumption~\ref{ass::fLip} holds, then we have
\[
\Phi \in \mathcal L( \; \mathcal P(X) \times \mathcal U \times \mathcal P(W), \; \mathcal P(X) \; ) \, ;
\]
that is, $\Phi$ is jointly Lipschitz continuous with respect to all three input arguments. Here, we recall that $\mathcal P(X)$ and $\mathcal P(W)$ are both metric spaces with respect to their associated Wasserstein distance functions $d_{\mathrm{W}}$, while $\mathcal U$ is equipped with the Hilbert-Sobolev distance function $d_{\mathcal L}$.
\end{lemma}

\pf First of all, as mentioned above, $\Phi$ is well defined by~\eqref{eq::PHI}, as Assumption~\ref{ass::fLip} ensures not only that $f$ is a Lipschitz continuous function but also that the image set of $f$ is contained in $X$---such that the image set of $\Phi$ is contained in $\mathcal P(X)$. Now, let $p,q \in \mathcal P(X)$ and $\omega,\xi \in \mathcal P(W)$ be given measures and $\mu,\nu \in \mathcal U$ given feedback laws. We set $p^+ = \Phi(p,\mu,\omega)$ and $q^+ = \Phi(q,\nu,\xi)$.
The definition of the Wasserstein metric implies that
\begin{eqnarray}
\begin{array}{l}
\displaystyle
d_{\mathrm{W}}(p^+,q^+) = \underset{\varphi \in \mathcal L_1(X,\mathbb R)}{\sup} \left| \int \varphi \, \mathrm{d}p^+ - \int \varphi \, \mathrm{d}q^+ \right| \notag \\[0.4cm]
\displaystyle
= \underset{\varphi \in \mathcal L_1(X,\mathbb R)}{\sup} \left| \int \int \varphi \circ f_{\mu} \, \mathrm{d}p \, \mathrm{d}\omega - \int \int \varphi \circ f_{\nu} \, \mathrm{d}q \, \mathrm{d}\xi \right|
\end{array}
\end{eqnarray}
Here, $\circ$ denotes the composition operator and, additionally, we have introduced the shorthand notation
\[
\forall \mu \in \mathcal U, \qquad f_{\mu}(x,w) \; \defeq \; f(x,\mu(x),w) \; .
\]
Because we assume that $\mathcal U \in \mathbb K( \mathcal L(X,U) )$ and because the particular definition of the Hilbert-Sobolev distance implies that all functions in $\mathcal U$ are uniformly Lipschitz continuous, the functions $f_{\mu}$ \mbox{are---by} \mbox{construction---uniformly} Lipschitz continuous over the compact set $\mathcal U$.
Let \mbox{$\gamma_1 < \infty$} denote the associated uniform Lipschitz constant. Since $\varphi$ is $1$-Lipschitz continuous, the functions of the form $\varphi \circ f_{\mu}$ are also Lipschitz continuous with uniform Lipschitz constant~\mbox{$1*\gamma_1 = \gamma_1$}. Thus, the triangle inequality yields the estimate 
\begin{eqnarray}
\begin{array}{l}
\displaystyle
\left| \int \int \varphi \circ f_{\mu} \, \mathrm{d}p \, \mathrm{d}\omega - \int \int \varphi \circ f_{\nu} \, \mathrm{d}q \, \mathrm{d}\xi \right| \\[0.4cm]
\displaystyle
\leq \; \left| \int \int \varphi \circ f_{\mu} \, \mathrm{d}p \, \mathrm{d}\omega - \int \int \varphi \circ f_{\mu} \, \mathrm{d}q \, \mathrm{d}\omega \right| \\[0.2cm]
\displaystyle
\hphantom{\leq \;} + \left| \int \int \varphi \circ f_{\mu} \, \mathrm{d}q \, \mathrm{d}\omega - \int \int \varphi \circ f_{\mu} \, \mathrm{d}q \, \mathrm{d}\xi \right| \\[0.2cm]
\displaystyle
\hphantom{\leq \;} + \left| \int \int \left[ \varphi \circ f_{\mu} - \varphi \circ f_{\nu} \right] \, \mathrm{d}q \, \mathrm{d}\xi\right| \\[0.4cm]
\displaystyle
\leq \; \gamma_1 \cdot d_{\mathrm{W}}(p,q) + \gamma_1 \cdot  d_{\mathrm{W}}(\omega,\xi) \\[0.25cm]
\displaystyle
\hphantom{\leq \;}  + \left| \int \int \left[ \varphi \circ f_{\mu} - \varphi \circ f_{\nu} \right] \, \mathrm{d}q \, \mathrm{d}\xi\right| \; ,
\end{array}
\end{eqnarray}
which holds uniformly for all $\varphi \in \mathcal L_1(X,\mathbb R)$ and all functions $\mu,\nu \in \mathcal U$. Additionally, since $q$ and $\xi$ are probability measures, the last integral term can be bounded as
\[
\left| \int \int \left[ \varphi \circ f_{\mu} - \varphi \circ f_{\nu} \right] \, \mathrm{d}q \, \mathrm{d}\xi\right| \leq \gamma_2 \cdot |X| \cdot |W| \cdot d_{\mathcal L}(\mu,\nu) \; ,
\]
where $\gamma_2$ denotes the Lipschitz constant of $f$ with respect to its second argument, $|X|$ the diameter of the compact set $X$ and $|W|$ the diameter of the set $W$. Finally, by substituting all the above inequalities we find that
\[
d_{\mathrm{W}}(p^+,q^+) \leq \gamma \left(
d_{\mathrm{W}}(p,q) + d_{\mathrm{W}}(\omega,\xi) + d_{\mathcal L}(\mu,\nu)
\right)
\]
with $\gamma = \max \{ \gamma_1 , \gamma_2 \cdot |X| \cdot |W| \}$. But this inequality implies that $\Phi$ is indeed Lipschitz continuous.\qed

\begin{remark}
The proof of Lemma~\ref{lem::Lipschitz} relies on the properties of Wasserstein (Kantorovich-Rubinstein) distances, which have originally been introduced independently by several authors including Kantorovich~\cite{Kantorovich2006} and Wasserstein~\cite{Vasershtein1969}, see also~\cite{Villani2005}. In order to explain why the Wasserstein metric is remarkably powerful in the context of control system analysis, let us briefly discuss what would have happened if we had used, for example, a total variation distance,
\[
d_{\mathrm{TV}}(p,q) \; \defeq \; \sup_{A \in \mathcal B(X)} | p(A)-q(A) | \; ,
\]
instead of the Wasserstein distance in order to define the metric space $\mathcal P(X)$. For this aim, we consider a scalar system with $f(x,u,w) = u$ and parametric feedback laws of the form $\mu_{\kappa}(x) = {\kappa}$, $\mathcal U = \{ \mu_{\kappa} \mid {\kappa} \in [-1,1] \}$. In this example, we have
\begin{eqnarray}
d_{\mathrm{TV}}(\Phi(p,\mu_{\kappa},\omega),\Phi(p,\mu_{{\kappa}'},\omega)) &=& d_{\mathrm{TV}}(\delta_{\kappa},\delta_{{\kappa}'})
\notag \\[0.16cm]
&=& \left\{
\begin{array}{ll}
0 \; &\text{if} \; {\kappa} = {\kappa}' \\
1 & \text{otherwise}
\end{array}
\right. \notag
\end{eqnarray}
implying that $\Phi$ is not Lipschitz continuous with respect to the parametric feedback law on $\mathcal U$ if we use the total variation distance to define our underlying metric space. In other words, the statement of Lemma~\ref{lem::Lipschitz} happens to be wrong in general, if we replace the Wasserstein distance with other distances such as the total variation distance.
\end{remark}

\subsection{Ambiguity Tubes}
\label{sec::AmbiguityTubes}
In this section, we generalize the considerations from the previous section by introducing ambiguity tubes. The motivation for considering such a general setting is twofold: firstly, in practice, one might not know the exact probability measure $\omega_k$ of the process noise $w_k$, but only have a set $\Omega \in \mathcal K(W)$ of possible probability measures. For instance, one might know a couple of lower order moments of $w_k$, such as the expected value and variance, while higher order moments are known less accurately or they could be even completely unknown. And secondly, in the context of nonlinear system analysis as well as in the context of high dimensional state spaces, a propagation of the exact state distribution can be difficult or impossible. In such cases, it may be easier to bound the true probability measure of the state by a so-called enclosure; that is, a set of probability measures that---in a suitable, yet to be defined sense---overestimates the actual probability measure of the state.

In order to prepare a mathematical definition of ambiguity tubes, we introduce an ambiguity transition map \mbox{$F: \mathcal K(X) \times \mathcal U \to \mathcal K(X)$}, which is defined as
\begin{eqnarray}
\label{eq::F}
F(P,\mu) \; \defeq \; \left\{ \; \Phi(p,\mu,\omega) \; \middle| \; p \in P, \; \omega \in \Omega \; \right\}
\end{eqnarray}
for all $P \in \mathcal K(X)$ and all $\mu \in \mathcal U$. Here, the ambiguity set \mbox{$\Omega \in \mathcal K(W)$} of possible disturbance probability measures is assumed to be given and constant. The following corollary is a direct consequence of Lemma~\ref{lem::Lipschitz}.

\begin{corollary}
\label{cor::Lipschitz}
Let Assumption~\ref{ass::fLip} hold.
The function $F$ is Lipschitz continuous,
\[
F \in \mathcal L( \; \mathcal K(X) \times \mathcal U, \; \mathcal K(X) \; ) \; ,
\]
recalling that $\mathcal K(X)$ is equipped with the Hausdorff-Wasserstein metric $d_\mathrm{H}$.
\end{corollary}

\pf First of all, the fact that the image sets of $F$ are compact follows from~\eqref{eq::F} and the fact that the function $\Phi$ is Lipschitz continuous. Next, $F$ directly inherits the Lipschitz continuity of~$\Phi$---including its Lipschitz constant~$\gamma$ that has been introduced in the proof of Lemma~\ref{lem::Lipschitz}---as $d_\mathrm{H}$ is the Hausdorff metric of $d_\mathrm{W}$.
\qed

In order to formalize the concept of ambiguity enclosures, the following definition is introduced.

\begin{definition}
\label{def::enclosure}
An ambiguity set $Q \in \mathcal K(X)$ is called an enclosure of an ambiguity set $P \in \mathcal K(X)$, denoted by $P \preceq Q$, if
\[
\sup_{\varphi \in \mathcal L_1(X,\mathbb R)} \left[ \max_{p \in P} \min_{q \in Q} \int \varphi \, \mathrm{d}(p-q) \right] \leq 0 \; .
\]
Two ambiguity sets $P,Q \in \mathcal K(X)$ are considered equivalent, denoted by \mbox{$P \simeq Q$}, if both $P \preceq Q$ and $Q \preceq P$.
\end{definition}

Notice that the above definition of the relation ``$\preceq$'' should not be mixed up with the set inclusion relation ``$\subseteq$'', as used for the definition of set enclosures in the field of set-based tube MPC and global optimization. The corresponding conceptual difference is illustrated by the following example.

\begin{example}
Let us consider the ambiguity sets
\[
P = \{ \delta_0 , \delta_1 \} \quad \text{and} \quad Q = \{ \delta_0, \delta_1, 0.5 \delta_0 + 0.5 \delta_1 \}
\]
of the compact set $X = [0,1] \subseteq \mathbb R$ recalling that $\delta_0$ and $\delta_1$ denote the Dirac measures at $0$ and $1$, respectively. In this example, the upper bounds of the integrals,
\begin{eqnarray}
\max_{p \in P} \int \varphi \; \mathrm{dp} &=& \max\{ \; \varphi(0), \; \varphi(1) \; \} \notag  \\[0.16cm]
\text{and} \qquad \max_{q \in Q} \int \varphi \; \mathrm{dq} &=& \max\{ \; \varphi(0), \; \varphi(1) \; \} \; , \notag
\end{eqnarray}
coincide for any Lipschitz continuous function $\varphi$. Similarly, the associated lower bounds
\begin{eqnarray}
\min_{p \in P} \int \varphi \; \mathrm{dp} &=& \min\{ \; \varphi(0), \; \varphi(1) \; \} \notag  \\[0.16cm]
\text{and} \qquad \min_{q \in Q} \int \varphi \; \mathrm{dq} &=& \min\{ \; \varphi(0), \; \varphi(1) \; \} \; , \notag
\end{eqnarray}
coincide, too. Consequently, in the sense of Definition~\ref{def::enclosure}, the ambiguity sets $P$ and $Q$ are equivalent, $Q \simeq P$. In particular, we have $Q \preceq P$ but we do not have $Q \subseteq P$. Thus, the relations $\preceq$ and $\subseteq$ are not the same. 
\end{example}

The following proposition ensures that the relation $\preceq$ defines a proper partial order on $\mathcal K(X)$ with respect to the equivalence relation $\simeq$. Moreover, topological compatibility with respect to our Hausdorff-Wasserstein metric setting is established.

\begin{proposition}
\label{prop::enclosure}
Let the enclosure relation $\preceq$ be defined as in Definition~\ref{def::enclosure}. Then, the following properties are satisfied for any $P,Q,T \in \mathcal K(X)$.

\begin{enumerate}

\addtolength{\itemsep}{2pt}

\item \textit{Reflexivity:} we have $P \preceq P$.

\item \textit{Anti-Symmetry:} if $P \preceq Q$ and $Q \preceq P$, then $P \simeq Q$.

\item \textit{Transitivity:} if $P \preceq Q$ and $Q \preceq T$, then $P \preceq T$.

\item \textit{Compactness:} The set $$\mathcal S = \{ (P,Q) \in \mathcal K(X) \times \mathcal K(X) \mid P \preceq Q \}$$ is compact; that is, $\mathcal S \in \mathbb K(\mathcal K(X) \times \mathcal K(X))$.
\end{enumerate}
\end{proposition}

\pf
Reflexivity, anti-symmetry with respect to the equivalence relation $\simeq$, and transitivity follow directly from Definition~\ref{def::enclosure}. Thus, our focus is on the last statement, which claims to establish compatibility of our definition of enclosures and the proposed Wasserstein-Hausdorff metric setting. Let $P_0,P_1,P_2, \ldots \in \mathcal K(X)$ and $Q_0,Q_1,Q_2, \ldots \in \mathcal K(X)$ be two convergent sequences with
$$P^\star \defeq \lim_{k \to \infty} P_k \qquad \text{and} \qquad Q^\star \defeq \lim_{k \to \infty} Q_k$$
and such that
$P_k \preceq Q_k$ for all $k \in \mathbb N$. Because all sets are compact the maximizers
\[
p_{k,\varphi}^\star \defeq \underset{p \in P_k}{\mathrm{argmax}} \int \varphi \, \mathrm{d}p \quad \text{and} \quad q_{k,\varphi}^\star \defeq \underset{q \in Q_k}{\mathrm{argmax}} \int \varphi \, \mathrm{d}q \; .
\]
exist for all $\varphi \in \mathcal L_1(X,\mathbb R)$. Next, since $\mathcal K(X)$ is a compact set of compact sets, we have not only $P^\star,Q^\star \in \mathcal K(X)$, but the triangle inequality for the Hausdorff-Wasserstein metric additionally yields that
\begin{eqnarray}
\widetilde p_{\varphi} \defeq \lim_{k \to \infty} p_{k,\varphi}^\star &=& \underset{p \in P^\star}{\mathrm{argmax}} \int \varphi \, \mathrm{d}p \notag \\[0.16cm] \text{and} \qquad \widetilde q_{\varphi} \defeq \lim_{k \to \infty} q_{k,\varphi}^\star &=& \underset{q \in Q^\star}{\mathrm{argmax}} \int \varphi \, \mathrm{d}q \; . \notag
\end{eqnarray}
A direct consequence of these equations is that we have
\begin{eqnarray}
& & \sup_{\varphi \in \mathcal L_1(X,\mathbb R)} \left[ \max_{p \in P^\star} \min_{q \in Q^\star} \int \varphi \, \mathrm{d}(p-q) \right] \notag \\[0.16cm]
& &\qquad = \sup_{\varphi \in \mathcal L_1(X,\mathbb R)} \{ \widetilde p_{\varphi} - \widetilde q_{\varphi} \} \notag \\[0.16cm]
& &\qquad = \sup_{\varphi \in \mathcal L_1(X,\mathbb R)} \lim_{k \to \infty} \{ \underbrace{p_{\varphi,k}^\star - q_{\varphi,k}^\star}_{\leq 0} \} \leq 0 \; , \notag
\end{eqnarray}
which shows that $P^\star \preceq Q^\star$. Notice that this means that if $(P_k,Q_k) \in S$ is a Cauchy sequence, then the limit point satisfies $(P^\star,Q^\star) \in S$; that is $S$ is closed. Because $S$ is bounded by construction, this also implies that $S$ is compact, \mbox{$S \in \mathbb K(\mathcal K(X) \times \mathcal K(X))$}.
\qed

\bigskip
\noindent
After this technical preparation, the following definition of ambiguity tubes is possible.

\begin{definition}
\label{def::tube}
The sequence $(P_0,P_1, \ldots, P_N) \in \mathcal K(X)^{N+1}$ is called an ambiguity tube of~\eqref{eq::system} on the discrete-time horizon $\{ 0,1, \ldots, N \}$, if there exist an associated  sequence of ancillary feedback controllers $\mu_0, \mu_1, \ldots, \mu_{N-1} \in \mathcal U$ such that
\[
\forall k \in \{ 0, 1, \ldots, N-1 \}, \qquad F(P_k,\mu_k) \preceq P_{k+1} \; .
\]
\end{definition}

Notice that Definition~\ref{def::tube} is inspired by methods from the field of set-theoretic methods in control, in particular, the concept of set-valued tubes, as used in the field of Tube MPC~\cite{Houska2019,Langson2004,Mayne2005}. In detail, the step from set-valued robust forward invariant tubes to ambiguity tubes is, however, not straightforward. For instance, the standard set inclusion relation ``$\subseteq$'' would be too strong for a practical definition of ambiguity tubes and is here replaced by the relation $\preceq$. This adaption of concepts to our measure based setting is needed, as the purpose of constructing tubes in traditional set propagation and ambiguity set propagation are different. As we will also discuss in the following sections, ambiguity tubes can be used to assess, analyze, and trade-off the risk of constraint violations with other performance measures rather than enforcing the conservative worst-case constraints that are implemented in traditional tube MPC.

\begin{remark}
An equivalent characterization of the relations in Definition~\ref{def::enclosure} can be obtained by borrowing notation from the field of convex optimization that is related to the concept of duality and support functions~\cite{Boyd2004,Rockafellar2005}. In order to explain this, we denote with $d_P: \mathcal L_1(X,\mathbb R) \to \mathbb R$ the support function of the ambiguity set $P$,
\[
\forall \varphi \in \mathcal L_1(X,\mathbb R), \qquad d_P(\varphi) \; \defeq \; \max_{p \in P} \int \varphi \, \mathrm{d}p \; .
\]
This notation is such that we have $d_P = d_Q$ if and only if $P \simeq Q$. Similarly, we have $d_P \leq d_Q$ if and only if $P \preceq Q$.
\end{remark}

\subsection{Proper Ambiguity Measures}
\label{sec::ambigiutyMeasures}
The goal of this section is to formalize certain concepts of modeling performance and risk in the space of ambiguity sets. For this aim, we introduce maps of the form
\[
\ell: \mathcal K(X) \to \mathbb R \; ,
\]
which assign real values to ambiguity sets. In this context, we propose to introduce the following regularity condition under which an ambiguity measure is considered ``proper''\footnote{The conditions in Definition~\ref{def::proper} are of an axiomatic nature that is inspired by similar axioms for coherent risk measures, as introduced in~\cite{Rockafellar2013}. One difference though is that we work with ambiguity sets rather than single probability measures. Moreover, Definition~\ref{def::proper} is, at least in this form, tailored to our proposed Wasserstein-Hausdorff metric setting in which Lipschitz continuity (not only closedness of image sets as required for regular risk measures~\cite{Rockafellar2013}) is needed for ensuring topological compatibility.}.

\begin{definition}
\label{def::proper}
We say that $\ell: \mathcal K(X) \to \mathbb R$ is a proper ambiguity measure, if $\ell$ is Lipschitz continuous, $\ell \in \mathcal L(\mathcal K(X),\mathbb R)$, linear with respect to weighted Minkowski sums; that is,
\[
\forall \theta \in [0,1], \quad \ell( \theta P \oplus (1-\theta) Q) = \theta \ell(P) + (1-\theta) \ell(Q) \; ,
\]
and monotonous; that is, $P \preceq Q$ implies $\ell(P) \leq \ell(Q)$ for $P,Q \in \mathcal K(X)$.
\end{definition}

The Lipschitz continuity requirement in the above definition ensures that proper ambiguity measures $\ell$ also satisfy the equation $\ell(P) = \ell(Q)$ for all $P,Q \in \mathcal K(X)$ with $P \simeq Q$. Together with the monotonicity relation, this implies that proper ambiguity measures are compatible with our definition of the relations ``$\preceq$'' and ``$\simeq$'' from Definition~\ref{def::enclosure}. In order to understand why practical performance and risk measures can be assumed to be proper without adding much of a restriction, it is helpful to have the following examples in mind.

\begin{example}
\label{ex::performance}
If $l \in \mathcal L(\mathbb R^n,\mathbb R)$ denotes a cost function, for example, the stage cost of a nominal MPC controller, the associated worst-case average performance
\[
\ell(P) \; \overset{\mathrm{def}}{=} \; \max_{p \in P} \int l \, \mathrm{d}p
\]
is well defined, where the maximizer exists for compact ambiguity sets $P \in \mathcal K(X)$. It is easy to check that $\ell$ is a proper ambiguity measure in the sense of Definition~\ref{def::proper}.
\end{example}

\begin{example}
\label{ex::risk}
If $\mathbb X \in \mathbb K(X)$ denotes a state constraint set, the associated maximum expected constraint violation at risk is given by
\[
\mathcal R(P) \; \overset{\mathrm{def}}{=} \; \max_{p \in P} \int \mathrm{dist}_1(x,\mathbb X) \, p(\mathrm{d}x)
\]
recalling that $\mathrm{dist}(x,\mathbb X) = \min_{z \in \mathbb X} \| x-z \|_1$ denotes the distance function with respect to the $1$-norm. Similar to the previous example, $\mathcal R$ turns out to be a proper ambiguity measure in the sense of Definition~\ref{def::proper}, which can here be interpreted as a risk measure. In fact, this risk measure is closely related to the so-called worst-case conditional value at risk, as introduced by T.R.~Rockafellar, see for example~\cite{Rockafellar2013}, which is by now accepted as one of the most practical and computationally tractable risk measures in engineering and management sciences.
\end{example}

\subsection{Ambiguity Tube MPC}
\label{sec::ambiguityMPC}
The focus of this section is on the formulation of ambiguity tube MPC problems of the form
\begin{eqnarray}
\mathcal V(y) \; \defeq &\underset{P,\mu}{\min}& \sum_{k=0}^{N-1} L(P_k,\mu_k) + M(P_N) \notag \\[0.16cm]
\label{eq::MPC}
&\mathrm{s.t.}& \left\{
\begin{array}{l}
\forall k \in \{ 0, 1, \ldots, N-1 \}, \\
F(P_k,\mu_k) \preceq P_{k+1} \\
P_k \in \mathcal K(X), \; \mu_k \in \mathcal U \\
\delta_y \in P_0 \; .
\end{array}
\right.
\end{eqnarray}
Here, the sequence of ambiguity sets $P = (P_0,P_1, \ldots, P_N)$ and the ancillary feedback laws $\mu = (\mu_0,\mu_1, \ldots, \mu_N)$ are optimization variables.
Moreover, $y$ denotes the current state \mbox{measurement---recalling} that $\delta_y$ denotes the associated Dirac \mbox{measure---while}
$$L: \mathcal K(X) \times \mathcal U \to \mathbb R \quad \text{and} \quad M: \mathcal K(X) \to \mathbb R$$ denote the stage and end cost terms. If $\mu_0^\star[y] \in \mathcal U$ denotes the parametric minimizer of~\eqref{eq::MPC}, then the actual MPC feedback law is given by
\begin{eqnarray}
\label{eq::MPCfeedbackLaw}
\mu_{\mathrm{MPC}}(y) \defeq \mu_0^\star[y](y) \; .
\end{eqnarray}
Notice that in this notation, the current time of the MPC controller is reset to $0$ after every iteration. The following theorem introduces a minimum requirement under which one could call~\eqref{eq::MPC} well-formulated.

\begin{theorem}
\label{thm::existence}
Let Assumption~\ref{ass::fLip} be satisfied and let $L$ and $M$ be continuous functions on the compact domains $\mathcal K(X) \times \mathcal U$ and $\mathcal K(X)$, respectively, then~\eqref{eq::MPC} admits a minimizer for any $y \in X$. Moreover, the function $\mathcal V$ is continuous on $X$.
\end{theorem}

\pf Because Assumption~\ref{ass::fLip} holds, we can combine the results of Corollary~\ref{cor::Lipschitz} with the fourth statement of Proposition~\ref{prop::enclosure} to conclude that the feasible set of~\eqref{eq::MPC} is non-empty and compact. Consequently, if $L$ and $M$ are continuous functions, Weierstrass' theorem yields the first statement of this theorem. Similarly, the second statement follows from a variant of Berge's theorem~\cite{Berge1963}; see, also~\cite[Thm~1.17]{Rockafellar2005}.
\qed

\bigskip
\noindent
Notice that the above theorem has been formulated under a rather weak requirement on the continuity of the functions $L$ and $M$; that is, without necessarily requiring that these functions are proper ambiguity measures. However, as we will discuss in the following sections, much stronger assumptions on the cost functions $L$ and $M$ are needed, if one is interested in analyzing the stability properties of the MPC controller~\eqref{eq::MPC}.

\begin{remark}
\label{rem::constraints}
Notice that the ambiguity tube MPC formulation~\eqref{eq::MPC} contains traditional tube MPC as well as existing stochastic MPC formulations as special cases that are obtained for the case that $\Omega$ is the set of all probability distributions with support set $W$ or, alternatively, a singleton, $\Omega = \{ \omega \}$. However, at this point, it should be mentioned that~\eqref{eq::MPC} is formulated in the understanding that state-constraints are taken into account by adding suitable risk measures to the objective, as explained by Example~\ref{ex::risk}. This notation \mbox{is---at} least from the viewpoint of stochastic \mbox{MPC---rather} natural, as one would in such a setting usually be interested in an MPC objective that allows one to tradeoff between the risk of violating a constraint and control performance. Nevertheless, for the sake of generality of the following analysis, it should be mentioned that if one is interested in enforcing explicit chance constraints on the probability of a constraint formulation, the corresponding MPC controllers can only be reformulated as a problem of the form~\eqref{eq::MPC}, if additional assumptions on the regularity\footnote{If we work with proper ambiguity measures in order to formulate constraints, this is clearly sufficient to ensure regularity.} and recursive feasibility of these constraints are \mbox{made---such} that they can be added to the stage cost in the form of $L_1$-penalties without altering the problem formulation. Notice that such regularity and recursive feasibility conditions have been discussed in all detail in~\cite{Rawlings2009} for min-max MPC and in~\cite{Kouvaritakis2016} for stochastic MPC.
\end{remark}

\section{Stability Analysis}
\label{sec::stability}
As mentioned in the introduction, the basic concepts for analyzing stability of Markovian systems have been established in~\cite{Bucy1965} and~\cite{Kushner1965} by using martingale theory. The goal of this section is to lay the foundation for applying this  theory to analyze the stochastic closed-loop stability properties of the ambiguity tube MPC controller~\eqref{eq::MPC} in the presence of uncertainties. For this aim, this section is divided into three parts: Section~\ref{sec::assumptions} concisely collects and elaborates on all assumptions that will be needed for this stability analysis, Section~\ref{sec::concave} establishes an important technical result regarding the concavity of MPC cost functions with respect to Minkowski addition of ambiguity sets, and Section~\ref{sec::supermartingale} uses this concavity property to construct a non-negative supermartingale, which finally leads to the powerful stability results for ambiguity tube MPC that are summarized in Theorems~\ref{thm::stability} and~\ref{thm::asympstability}.

\subsection{Conditions on the Stage and Terminal Cost Function}
\label{sec::assumptions}
Throughout the following stability analysis, two main assumptions on the stage and end cost function of the MPC controller~\eqref{eq::MPC} are needed, as introduced below.

\begin{assumption}
\label{ass::proper}
The functions $L(\cdot,\mu)$ and $M$ are for any given $\mu \in \mathcal U$ proper ambiguity measures. Moreover, we assume that $L$ is continuous on $\mathcal K(X) \times \mathcal U$.
\end{assumption}

\noindent
Notice that Examples~\ref{ex::performance} and~\ref{ex::risk} discuss how to formulate practical risk and performance measures in such a way that the above assumption holds. From here on, a separate assumption on the continuity of $M$ (as in Theorem~\ref{thm::existence}) is not needed anymore, as proper ambiguity measures are Lipschitz continuous functions and, consequently, also continuous. Assumption~\ref{ass::proper} does, however, add the explicit requirement that $L$ is continuous, as this function depends in general also on the feedback \mbox{law~$\mu$---in} this way, we make sure that the conditions of Theorem~\ref{thm::existence} are satisfied whenever Assumptions~\ref{ass::fLip} and~\ref{ass::proper} are satisfied. The following assumption formulates an additional condition under which we consider the terminal cost function $M$ admissible.

\begin{assumption}
\label{ass::terminal}
The functions $L$ and $M$ are non-negative and satisfy the terminal descent condition
\[
\forall P \in \mathcal K(X), \; \exists \mu \in \mathcal U: \; \; L(P,\mu) + M( F(P,\mu) ) \leq M(P) \, .
\]
\end{assumption}

\noindent
The above assumption can be interpreted as a Lyapunov decent \mbox{condition---similar} conditions are used for constructing terminal cost functions for classical certainty-equivalent as well as tube based MPC controllers~\cite{Chen1998,Gruene2009,Rawlings2009,Villanueva2020b}. The question how to construct functions $L$ and $M$ that simultaneously satisfy Assumptions~\ref{ass::proper} and~\ref{ass::terminal} in practice will be discussed in Section~\ref{sec::terminalCost}.

\begin{remark}
Assumptions~\ref{ass::proper} and~\ref{ass::terminal} together imply that $M$ must be a Lipschitz continuous Control Lyapunov Function (CLF). In the general context of traditional nonlinear system analysis, conditions under which such Lipschitz continuous CLFs exist have been analyzed by various authors~\cite{Clarke1998,Ledyaev1999}. However, if one considers nonlinear MPC problems with explicit state constraints (see also Remark~\ref{rem::constraints}), it is possible to construct \mbox{systems---for} example, based on Artstein's \mbox{circles---that} are asymptotically stabilizable yet fail to admit a continuous CLF~\cite{Grimm2004}. It is, however, also pointed out in~\cite{Grimm2004} that systems that only admit discontinuous CLFs often lead to non-robust MPC controllers. Therefore, in the context of robust MPC design, the motivation for working with Assumptions~\ref{ass::proper} and~\ref{ass::terminal} is to exclude such pathological non-robustly stabilizable systems. Notice that more general regularity assumptions, under which a robust control design is possible, are beyond the scope of this paper.
\end{remark}

\subsection{On Concave Cost Functions}
\label{sec::concave}
After summarizing all main assumptions of this section, we can now focus on the properties of certain cost functions that will later be used to construct a supermartingale for the proposed ambiguity tube MPC controller. For this aim, we first introduce the auxiliary function
\begin{eqnarray}
\mathcal{J}_\mu(Q) \; \defeq & \underset{P}{\min} & \sum_{k=0}^{N-1} L(P_k,\mu_k) + M(P_N) \notag \\[0.16cm]
\label{eq::AuxCost}
&\mathrm{s.t.}& \left\{
\begin{array}{l}
\forall k \in \{ 0, 1, \ldots, N-1 \}, \\
F(P_k,\mu_k) \preceq P_{k+1} \\
P_k \in \mathcal K(X) \\
Q \preceq P_0, \; P_N = P^\star  \; ,
\end{array}
\right.
\end{eqnarray}
which is formally defined for all $Q \in \mathcal K(X)$ and all ancillary feedback laws $\mu \in \mathcal U^N$.

\begin{lemma}
\label{lem::linJaux}
If Assumptions~\ref{ass::fLip} and~\ref{ass::proper} are satisfied, then $\mathcal J_\mu$ is for any given $\mu \in \mathcal U^N$ a proper ambiguity measure.
\end{lemma}

\pf
In the following, we may assume that $\mu \in \mathcal U^N$ is constant and given. Our proof is divided into two parts, where the first part focuses on establishing a linearity property of the function $F$. The second part of the proof builds upon the first part in order to further analyze the properties of the function $\mathcal J_{\mu}$.

\bigskip
\noindent
\textit{PART I:} Recall that the function $\Phi$, which has been defined in~\eqref{eq::PHI}, is---by construction---linear in its first argument. Consequently, we have
\[
\theta \Phi(p,\nu,\omega) + (1-\theta) \Phi(q,\nu,\omega) = \Phi(\theta p + (1-\theta) q,\nu,\omega)
\]
for all $p,q \in \mathcal P(X)$ and all $\theta \in [0,1]$, for any given \mbox{$\omega \in \mathcal P(W)$} and \mbox{$\nu \in \mathcal U$}. But this property of $\Phi$ implies directly that the function $F$ satisfies
\[
\theta F(P,\nu) \oplus (1-\theta) F(Q,\nu) = F(\theta P \oplus (1-\theta) Q,\nu) \; ,
\]
for all $P,Q \in \mathcal K(X)$, all $\nu\in \mathcal U$, and all $\theta \in [0,1]$, which follows from the definition of $F$ in~\eqref{eq::F}. 

\bigskip
\noindent
\textit{PART II:} Let $\mathfrak P_0[Q], \ldots, \mathfrak P_N[Q] \in \mathcal K(X)$ denote the solution of the recursion
\begin{eqnarray}
\mathfrak{P}_0[Q] &=& Q \; , \notag \\[0.16cm]
\mathfrak P_{k+1}[Q] &=& F( \mathfrak P_k[Q], \mu_k ) \notag
\end{eqnarray}
for $k \in \{ 0,1, \ldots, N-1\}$ recalling that $\mu$ is given. Corollary~\ref{cor::Lipschitz} ensures that the transition map $F$ is Lipschitz continuous such that the above recursion generates indeed compact ambiguity sets for any compact input set $Q \in \mathcal K(X)$, such that the sequence $\mathfrak P_0, \ldots, \mathfrak P_N$ is well-defined. Due to the linearity of $F$ with respect to its first argument (see Part~I), it follows by induction over $k$ that
\begin{eqnarray}
\label{eq::PP1}
\mathfrak P_k[\theta Q \oplus (1-\theta)Q'] = \theta \mathfrak P_k[Q] \oplus (1-\theta) \mathfrak P_k[Q']
\end{eqnarray}
for all $Q,Q' \in \mathcal K(X)$ and all $\theta \in [0,1]$. Next, because both the ambiguity measures $L$ and $M$ \mbox{are---due} to \mbox{Assumption~\ref{ass::proper}---monotonous} and Lipschitz continuous, we have
\begin{eqnarray}
\label{eq::JJ1}
\mathcal J_{\mu}(Q) = \sum_{k=0}^{N-1} L(\mathfrak P_k[Q],\mu_k) + M(\mathfrak P_N[Q])
\end{eqnarray}
for any $Q \in \mathcal K(X)$.

\smallskip
\noindent
Finally, it remains to combine the results~\eqref{eq::PP1} and~\eqref{eq::JJ1} with our assumption that $L$ and $M$ are proper ambiguity measures, which implies that $\mathcal J_{\mu}$ is a proper ambiguity measure, too.
\qed

\bigskip
\noindent
Notice that the auxiliary function $\mathcal J_\mu$ depends on the feedback law $\mu$, which is optimized in the context of MPC. Consequently, we are in the following not directly interested in this auxiliary function, but rather in the actual cost-to-go function
\begin{eqnarray}
\label{eq::CostToGo}
J(Q) \; \defeq \; \min_{\mu \in \mathcal U^N} \; \mathcal J_{\mu}(Q) \; ,
\end{eqnarray}
which is defined for all $Q \in \mathcal K(X)$, too. Clearly, the function $J$ is closely related to the value function $\mathcal V$ of the ambiguity tube MPC controller~\eqref{eq::MPC}, as we have
\begin{eqnarray}
\label{eq::VJ}
\forall y \in X, \qquad \mathcal V(y) = J( \{ \delta_y \} ) \; .
\end{eqnarray}
This equation follows directly by comparing the definition of $\mathcal V$ in~\eqref{eq::MPC} with the definitions in~\eqref{eq::AuxCost} and~\eqref{eq::CostToGo}.
The following corollary summarizes an important consequence of Lemma~\ref{lem::linJaux}.

\begin{corollary}
\label{cor::concave}
Let Assumptions~\ref{ass::fLip} and~\ref{ass::proper} be satisfied. The cost-to-go function $J$ is concave with respect to weighted Minkowski addition; that is, we have
\begin{eqnarray}
\label{eq::concave}
J( \theta Q \oplus (1-\theta)Q' ) \; \geq \; \theta J( Q ) + (1-\theta) J( Q' )
\end{eqnarray}
for all $Q,Q' \in \mathcal K(X)$ and all $\theta \in [0,1]$. Moreover, $J$ is monotonous, $Q \preceq Q'$ implies $J(Q) \leq J(Q')$. 
\end{corollary}

\pf
The key idea for establishing the first statement of this corollary is to use the linearity of the auxiliary function $\mathcal J_{\mu}$ (Lemma~\ref{lem::linJaux}). This yields the inequality
\begin{eqnarray}
\begin{array}{l}
J( \theta Q \oplus (1-\theta)Q' ) \\[0.25cm]
\begin{array}{cl}
=& \underset{\mu \in \mathcal U^N}{\min} \mathcal J_{\mu}( \theta Q \oplus (1-\theta)Q' ) \notag \\[0.3cm]
=& \underset{\mu \in \mathcal U^N}{\min} \left\{ 
\theta \mathcal J_\mu(Q) + (1-\theta)\mathcal J_{\mu}(Q')
\right\} \notag \\[0.3cm]
\geq& 
\theta \left( \underset{\mu \in \mathcal U^N}{\min} \, \mathcal J_{\mu}(Q) \right) + (1-\theta)\left( \underset{\mu' \in \mathcal U^N}{\min} \, \mathcal J_{\mu'}(Q') \right) \notag \\[0.5cm]
=& \theta J( Q ) + (1-\theta) J( Q' )
\end{array}
\end{array}
\end{eqnarray}
for all $Q,Q' \in \mathcal K(X)$ and all $\theta \in [0,1]$. This corresponds to the first statement of the lemma. The second statement follows from the fact that any minimizer of the optimization problem
\begin{eqnarray}
\mathcal{J}_\mu(Q') \; = & \underset{P,\mu}{\min} & \sum_{k=0}^{N-1} L(P_k,\mu_k) + M(P_N) \notag \\[0.16cm]
\label{eq::AuxCost22}
&\mathrm{s.t.}& \left\{
\begin{array}{l}
\forall k \in \{ 0, 1, \ldots, N-1 \}, \\
F(P_k,\mu_k) \preceq P_{k+1} \\
P_k \in \mathcal K(X) \; , \; \mu_k \in \mathcal U \\
Q' \preceq P_0 \; ,
\end{array}
\right.
\end{eqnarray}
is a feasible point of the corresponding optimization problem that is obtained when replacing the ambiguity set $Q' \in \mathcal K(X)$ with an ambiguity set $Q$ that satisfies $Q \preceq Q'$, which then implies monotonicity; that is, \mbox{$J(Q) \leq J(Q')$}.
\qed

\bigskip
\noindent
Corollary~\ref{cor::concave} is central to establishing stability properties of stochastic MPC or, in the context of this paper, more general ambiguity tube MPC schemes. The following example helps to understand why this concavity statement is important.

\begin{example}
Let us consider the example that $Q = \{ \delta_a \}$ and $Q' = \{ \delta_b \}$ for two points $a,b \in X$. In words, this means that $J(Q)$ is the cost that is associated with knowing that we are currently at the point $a$. Similarly, $J(Q')$ can be interpreted as the cost that is associated with knowing that we are currently at the point $b$. Now, if we set $\theta = \frac{1}{2}$, the corresponding ambiguity set
\[
\theta Q \oplus (1-\theta)Q' = \left\{ \frac{1}{2} \delta_a + \frac{1}{2} \delta_b \right\}
\]
can be associated with the situation that we don't know whether we are at the point $a$ or at the point $b$, as both events could happen with probability $\frac{1}{2}$. Thus, in this example, the first statement of Corollary~\ref{cor::concave} is saying that the cost that is associated with not knowing whether we are at $a$ \mbox{or $b$---that} is, \mbox{$J(\theta Q \oplus (1-\theta)Q')$---is} larger or equal than the expected cost that is obtained when planning to first measure whether we are at $a$ or $b$ and then evaluating the cost function.
\end{example}

\smallskip
\noindent
The following section exploits the conceptual idea from the above example and Corollary~\ref{cor::concave} in order to construct a non-negative supermartingale for ambiguity tube MPC.

\subsection{Supermartingales for Ambiguity Tube MPC}
\label{sec::supermartingale}
The goal of this section is to establish conditions under which the value function $\mathcal V$ is a supermartingale along the trajectories of the closed system that is associated with the ambiguity tube MPC feedback law $\mu_{\mathrm{MPC}}$, as defined in~\eqref{eq::MPCfeedbackLaw}. Let us first recall that
\begin{eqnarray}
\label{eq::ClosedLoop}
\forall k \in \mathbb N, \qquad y_{k+1} = f( y_k, \mu_{\mathrm{MPC}}(y_k), w_k  )
\end{eqnarray}
denotes the closed-loop stochastic process that is associated with the MPC controller~\eqref{eq::MPC}. Here, we additionally recall that $w_0, w_1, \ldots : W \to W$, $w_k(w)=w$, are independent $B(W)$-measurable random variables in the probability spaces $(W,\mathcal B(W),\omega_k)$ that depend on the sequence of measures $\omega_0,\omega_1, \ldots \in \Omega$. This implies that the states $y_0,y_1, \ldots$ are random variables, too. In the following, we use the notation $\mathfrak S_k = \sigma( y_0, y_1, \ldots, y_k )$ to denote the minimal \mbox{$\sigma$-field} of the sequence $y_0,y_1, \ldots, y_k$, such that
\[
\forall k \in \mathbb N, \qquad \mathfrak S_k \subseteq \mathfrak S_{k+1}
\]
is a filtration of $\sigma$-algebras. Moreover, we use the standard notation~\cite{Taylor1996}
\begin{eqnarray}
\mathbb E \{ \mathcal V(y_{k+1}) \mid \mathfrak S_k \} &\defeq& \int \mathcal V(f( y_k, \mu_{\mathrm{MPC}}(y_k), \cdot  )) \, \mathrm{d} \omega_k \notag
\end{eqnarray}
to denote the conditional expectation of $\mathcal V(y_{k+1})$ given~$\mathfrak S_k$. Notice that this definition depends on the sequence of probability measures $\omega_0, \omega_1, \ldots \in \Omega$. Clearly, if Assumptions~\ref{ass::fLip}, and~\ref{ass::proper} hold, Theorem~\ref{thm::existence} ensures that $\mathcal V$ is continuous. Since Assumption~\ref{ass::fLip} also ensures that $f$ is Lipschitz continuous, the integrand in above expression is trivially $\mathcal B(W)$-measurable and, consequently, the conditional expectation of $\mathcal V(y_{k+1})$ given $\mathfrak S_k$ is well-defined.

\begin{theorem}
\label{thm::martingale}
Let Assumptions~\ref{ass::fLip},~\ref{ass::proper}, and~\ref{ass::terminal} be satisfied. Then the function $\mathcal V$ is a supermartingale along the trajectories of the closed-loop system~\eqref{eq::ClosedLoop}; that is, we have
\[
\forall k \in \mathbb N, \qquad \mathbb E\{ \mathcal V (y_{k+1}) \mid \mathfrak S_k \} \leq \mathcal V(y_k)
\]
independent of the choice of $\omega_0,\omega_1 \ldots \in \Omega$.
\end{theorem}

\pf
Throughout this proof we use the notation $P_0^\star(y), P_1^\star(y), \ldots \in \mathcal K(X)$ to denote the parametric minimizers of~\eqref{eq::MPC} such that the inequality
\begin{eqnarray}
\label{eq::Vup}
\mathbb E \{ \mathcal V (y_{k+1}) \mid \mathfrak S_k \} &\leq& \max_{p \in P_1^\star(y_n)} \, \int \mathcal V \, \mathrm{d}p
\end{eqnarray}
holds---by definition of~$P_1^\star(y_n)$---for any choice of the sequence $\omega_0,\omega_1 \ldots \in \Omega$. Next, since Assumptions~\ref{ass::fLip} and~\ref{ass::proper} are satisfied, we can use the first statement of Corollary~\ref{cor::concave} to establish the inequality
\begin{eqnarray}
\int \mathcal V \, \mathrm{d}p &\overset{\eqref{eq::VJ}}{=}& \int J \left(  \{ \delta_y \} \right) \, p(\mathrm{d}y ) \notag \\[0.16cm]
\label{eq::jensen}
&\overset{\eqref{eq::concave}}{\leq}& J\left( \, \left\{ \int \delta_y \, p(\mathrm{d}y) \right\} \, \right) = J(\{ p \} ) \; ,
\end{eqnarray}
which holds for any probability measure $p \in \mathcal P(X)$. Moreover, we know from the second statement of Corollary~\ref{cor::concave} that $J$ is monotonous, which implies that the implication chain
\begin{eqnarray}
\label{eq::implic}
p \in P \quad \Longrightarrow \quad \{ p \} \preceq P \quad \Longrightarrow \quad J(\{ p \}) \leq J(P)
\end{eqnarray}
holds for any $P \in \mathcal K(X)$. In order to briefly summarize our intermediate results so far, we combine~\eqref{eq::Vup},~\eqref{eq::jensen}, and~\eqref{eq::implic}, which yields the inequality
\begin{eqnarray}
\label{eq::Jup}
\mathbb E \{ \mathcal V (y_{k+1}) \mid \mathfrak S_k \} &\leq& J( P_1^\star(y_k) ) \; .
\end{eqnarray}
Next, in analogy to the construction of Lyapunov functions for traditional MPC controllers~\cite{Rawlings2009}, we can use that Assumption~\ref{ass::terminal} implies that the cost-to-go function $J$ descends along the iterates $P_i^\star(y_k)$, which means that
\begin{eqnarray}
\label{eq::Jdecent}
J( P_1^\star(y_k)) &\leq& J( P_0^\star(y_k)) \; = \; \mathcal V(y_k) \; .
\end{eqnarray}
But then~\eqref{eq::Jup} and~\eqref{eq::Jdecent} imply that we also have
\[
\mathbb E \{ \mathcal V (y_{k+1}) \mid \mathfrak S_k \} \; \leq \; \mathcal V(y_k) \; .
\]
The latter inequality does not depend on the choice of the sequence $\omega_0,\omega_1, \ldots \in \Omega$ and, consequently, corresponds to the statement of the theorem.
\qed

\subsection{Stability of Ambiguity Tube MPC}
As established in Theorem~\ref{thm::existence}, Assumptions~\ref{ass::fLip} and~\ref{ass::proper} are sufficient to ensure that $\mathcal V$ is a continuous function. Consequently,
\begin{eqnarray}
\label{eq::Ydef}
\mathcal V^\star \; \defeq \; \min_{y \in X} \, \mathcal V(y) \quad \text{and} \quad Y^\star \; \defeq \; \underset{y \in X}{\text{argmin}} \;\mathcal V(y) 
\end{eqnarray}
are well-defined and the set $Y^\star$ is compact, $Y^\star \in \mathbb K(X)$, and non-empty (see Theorem~\ref{thm::existence}). Moreover, we may assume,  without loss of generality, that $\mathcal V^\star = 0$, as adding constant offsets to $\mathcal V$ does not affect the result of Theorem~\ref{thm::martingale}. In the following, we introduce the notation
\begin{eqnarray}
\mathcal N_{ \varepsilon }(Y^\star) &\defeq& \left\{ \; x \in X \; \middle| \; \mathrm{dist}_2(x,Y^\star) < \varepsilon \; \right\} \notag
\end{eqnarray}
to denote an $\varepsilon$-neighborhood of $Y^\star$. The definition below introduces a (rather standard) notion of stability of the closed-loop system~\eqref{eq::ClosedLoop}. Here, it is helpful to recall that the probability measures $p_k$ of the sequence $y_k$ are given by the recursion
\[
p_{k+1} = \Phi(p_k, \mu_{\mathrm{MPC}},\omega_k) \; ,
\]
which, in turn, depends on the sequence $\omega_0,\omega_1, \ldots \in \Omega$ and on the given initial state measurement $y_0$, as we set $p_0 = \delta_{y_0}$.

\begin{definition}
The closed-loop system~\eqref{eq::ClosedLoop} is called robustly stable with respect to the set $Y^\star$, if there exists for every $\varepsilon > 0$ a $\delta > 0$ such that for any $y_0 \in \mathcal N_{\delta}(Y^\star)$ we have
\[
\forall k \in \mathbb N, \quad p_k \left( \mathcal N_{\varepsilon}(Y^\star) \right) \; > \; 1-\varepsilon \, ,
\]
independent of the choice of the probability measures $\omega_0, \omega_1, \ldots \in \Omega$. If we additionally have that
\[
\lim_{k \to \infty} \; p_k( Y^\star ) = 1 \; ,
\]
independent of the choice of $\omega_0, \omega_1, \ldots \in \Omega$, we say that the closed-loop system is robustly asymptotically stable with respect to $Y^\star$.
\end{definition}

\noindent
The following theorem summarizes the main result of this section.

\begin{theorem}
\label{thm::stability}
If Assumptions~\ref{ass::fLip},~\ref{ass::proper}, and~\ref{ass::terminal} hold, then~\eqref{eq::ClosedLoop} is robustly stable with respect to $Y^\star$.
\end{theorem}

\pf
The assumptions of this theorem ensure that $\mathcal V$ is continuous (Theorem~\ref{thm::existence}) and a supermartingale along the trajectories of~\eqref{eq::ClosedLoop}, independent of the choice of the probability measures $\omega_0, \omega_1, \ldots \in \Omega$ (Theorem~\ref{thm::martingale}). Moreover, since we work with a compact support, all random variables are essentially bounded. Consequently, we can apply Bucy's supermartingale stability theorem~\cite[Thm.~1]{Bucy1965} to conclude that~\eqref{eq::ClosedLoop} is robustly stable with respect to $Y^\star$.
\qed

\begin{remark}
\label{rem::bucy}
The above proof is based on the historical result of Bucy's original article on positive supermartingales, \mbox{who---at} the time of publishing his original \mbox{article---formally} only established stability of Markov processes with respect to an isolated equilibrium; that is, for the case that the set $Y^\star$ is a singleton. However, the proof of Theorem~1 in~\cite{Bucy1965} generalizes trivially to the version that is needed in the above proof after replacing the distance between the states of the Markov system and the equilibrium point by the corresponding distance of these iterates to the set $Y^\star$. By now, this and other generalization of the supermartingale based  stability theorems by Bucy and Kushner for Markov processes are, of course, well-known and can in very similar versions also  (and besides many others) be found in~\cite{Feller1971,Kushner1965,Kushner2014,Taylor1996}.
\end{remark}

For the case that we are not only interested in robust stability but also robust asymptotic stability, the above theorem can be extended after introducing a slightly stronger regularity requirement on the function $L$, which is introduced below.

\begin{definition}
We say that the function $L$ is positive definite with respect to $Y^\star$ if $L(P,\mu) > 0$ for all $P$ with
$$\min_{p \in P} \; p(Y^\star) < 1$$
and all $\mu \in \mathcal U$.
\end{definition}

The corresponding stronger version of Theorem~\ref{thm::stability} can now be formulated as follows.

\begin{theorem}
\label{thm::asympstability}
Let Assumptions~\ref{ass::fLip},~\ref{ass::proper}, and~\ref{ass::terminal} hold. If $L$ is positive definite with respect to $Y^\star$, then~\eqref{eq::ClosedLoop} is robustly asymptotically stable with respect to $Y^\star$.
\end{theorem}

\pf
The statement of this theorem is similar to the statement of Theorem~\ref{thm::stability}, but we need to work with a slightly tighter version of the supermartingale inequality from Theorem~\ref{thm::martingale}. For this aim, we use that Assumption~\ref{ass::terminal} implies that
\begin{eqnarray}
\label{eq::Jdecent2}
L(P_0^\star(y_k),\mu^\star[y_k]) + J( P_1^\star(y_k)) &\leq& J( P_0^\star(y_k)) \; .
\end{eqnarray}
Consequently, the inequality~\eqref{eq::Jdecent} can be replaced by its tighter version,
\[
\min_{p \in P_0^\star(y_k)} \; p(Y^\star) < 1 \quad \Longrightarrow \quad J( P_1^\star(y_k)) < \mathcal V(y_k) \; .
\]
Thus, the corresponding argument in the proof of Theorem~\ref{thm::martingale} can be modified finding that we also have
\[
\min_{p \in P_0^\star(y_k)} \; p(Y^\star) < 1 \quad \Longrightarrow \quad \mathbb E\{ \mathcal V (y_{k+1}) \mid \mathfrak S_k \} < \mathcal V(y_k)
\]
for all $k \in \mathbb N$ and independent of the choice of the sequence $\omega_0,\omega_1 \ldots \in \Omega$. But this means that $\mathcal V$ is a strict non-negative supermartingale and we can apply the standard result from~\cite[Thm.~2]{Bucy1965} (of course, again after replacing Bucy's outdated definition of equilibrium points with our definition of the set $Y^\star$---see Remark~\ref{rem::bucy}) to establish the statement of this theorem.
\qed

\section{Practical Implementation of Ambiguity Tube MPC for Linear Systems}
\label{sec::implementation}
In order to illustrate and discuss the above theoretical results, this section develops a practical framework for reformulating a class of ambiguity tube MPC controllers for linear systems as convex optimization problems that can then be solved with existing MPC software. In particular, Section~\ref{sec::terminalCost} focuses on the question how to construct stabilizing terminal costs for stochastic and ambiguity tube MPC.

\subsection{Linear Stochastic Systems}
Let us consider linear stochastic discrete-time systems with (projected) linear right-hand function
\[
f = \mathrm{proj}_X \circ \hat f \qquad \text{with} \qquad \hat f(x,u,w) = A x + B u + w \; ,
\]
where $A \in \mathbb R^{n_x \times n_x}$ and $B \in \mathbb R^{n_x \times n_u}$ are given matrices. Here, we assume that an asymptotically stabilizing linear feedback gain $K \in \mathbb R^{n_u \times n_x}$ for $(A,B)$ exists; that is,
such that all eigenvalues of $A+BK$ are in the open unit disc. For simplicity of presentation, we focus on parametric ancillary feedback controllers of the form
\[
\mu[u_\mathrm{c}](x) = u_\mathrm{c} + Kx \; , \quad \text{where} \quad \mathcal U = \left\{ \; \mu[u_\mathrm{c}] \; \mid \; u_\mathrm{c} \in U_\mathrm{c} \; \right\}
\]
denotes the associated compact set of representable ancillary control laws for the pre-computed control gain $K$. In this context, $U_\mathrm{c} \in \mathbb K(\mathbb R^{n_u})$ denotes a control constraint that is associated with the central control offset~$u_\mathrm{c}$. Similarly, we introduce the notation
\[
\omega[ w_\mathrm{c} ](X') = \overline \omega( \{ w_\mathrm{c} \} \oplus X') \; \; \text{and} \; \; \Omega = \{ \omega[ w_\mathrm{c} ] \mid w_\mathrm{c} \in W_\mathrm{c} \}
\]
for all $X' \in \mathcal B(X)$,
assuming that $\overline \omega \in \mathcal P(\overline W)$, the central set $W_\mathrm{c} \in \mathbb K( \mathbb R^{n_x} )$, and the domain $\overline W \in \mathbb K( \mathbb R^{n_x} )$ are given. Finally, the corresponding domain of the uncertainty sequence is denoted by
\[
W \; \defeq \; W_\mathrm{c} \oplus \overline W \; .
\]
At this point, there are two remarks in order.

\begin{remark}
Because $K$ is a pre-stabilizing feedback, one can construct the compact domain $X \in \mathbb K( \mathbb R^{n_x} )$ such that
\[
X \supseteq (A+BK) X \oplus \left[ B U_\mathrm{c} \oplus W \right] \; .
\]
This construction is such that the closed-loop uncertainty measure propagation operator $\Phi$ is unaffected by the projection \mbox{onto $X$---it} would have been the same, if we would have set $f = \hat f$. This illustrates how Assumption~\ref{ass::fLip} \mbox{can---at} least for asymptotically stabilizable linear \mbox{systems---formally} be satisfied by a simple projection onto an invariant set, without altering the original physical problem formulation.
\end{remark}

\begin{remark}
\label{rem::nonlinear}
Notice that the ambiguity set $\Omega$ in the above system model can be used to overestimate nonlinear terms. For example, if we have a system of the form 
\[
x^+ = Ax + Bu + g(x) + \overline{w} \; ,
\]
where $g(x)$ is bounded on $X$, such that $g(x) \in W_\mathrm{c}$ while $\overline w$ is a random variable with probability measure $\overline \omega$, then the ambiguity set $\Omega$ is such that the probability measure $\omega_g$ of the random variable $w = g(x)+\overline w$ satisfies $\{ \omega_g \} \preceq \Omega$. This example can be used as a starting point to develop computationally tractable ambiguity tube MPC formulations for nonlinear systems, although a discussion of less conservative nonlinearity bounders, as developed for set propagation in~\cite{Villanueva2017}, are beyond the scope of this paper. 
\end{remark}

\subsection{Construction of Stage and Terminal Costs}
\label{sec::terminalCost}
In order to discuss how to design stage and terminal costs, $L$ and $M$, which satisfy the requirements from Assumptions~\ref{ass::proper} and~\ref{ass::terminal}, this section focuses on the case that the stage cost has the form
\begin{eqnarray}
\label{eq::LL}
L(P,\mu) = \max_{p \in P} \int \Theta(x,u_\mathrm{c}) \, p(\mathrm{d}x) \; ,
\end{eqnarray}
where $\Theta \in \mathcal L( X \times U_\mathrm{c}, \mathbb R_+ )$ is a non-negative and Lipschitz continuous control performance function. In the following, we assume that we have $0 \in U_\mathrm{c}$ as well $\Theta(y^\star,0) = 0$ for all $y^\star \in Y^\star$, where $Y^\star$ denotes the limit set of the considered linear stochastic system that is obtained for the offset-free ancillary control law; that is,
\[
Y^\star \; = \; \bigoplus_{k=0}^{\infty} (A+BK)^k \, W \; .
\]
Notice that this set can be computed by using standard methods from the field of set based computing~\cite{Blanchini2009,Houska2019}. 

\begin{example}
\label{ex::theta}
Let us assume that $\mathbb X \in \mathbb K(\mathbb R^{n_x})$ is a given state constraint with $Y^\star \subseteq \mathbb X$. In this case, the risk and performance measure
\[
\Theta(x,u) = \| u \|_2^2 + \mathrm{dist}_{2}(x,Y^\star)^2 + \tau \cdot \mathrm{dist}_{1}(x,\mathbb X)
\]
is Lipschitz continuous. It can be used to model a trade-off between the least-squares control performance term
\[
\| u \|_2^2 + \mathrm{dist}_{2}(x,Y^\star)^2
\]
that penalizes control offsets and the distance of the state to the target region $Y^\star$, and the constraint violation term $\tau \cdot \mathrm{dist}_{1}(x,\mathbb X)$ that is $0$ if $x$ satisfies the constraint. Here, $\tau > 0$ is a tuning parameter that can be used to adjust how risk-averse the controller is; see also Examples~\ref{ex::performance} and~\ref{ex::risk}.
\end{example}

Now, the key idea for constructing the terminal cost $M$ is to first construct a non-negative function $\Pi \in \mathcal L(X,\mathbb R_+)$ that satisfies the ancillary Lyapunov descent condition
\begin{eqnarray}
\label{eq::LMd}
\forall x \in X, \quad
\Theta(x,0) + \underset{w \in W}{\max} \; \Pi( f(x,Kx,w) ) \; \leq \; \Pi(x) \; .
\end{eqnarray}
Notice that such a function $\Pi$ exists, as we assume that the closed-loop system matrix $(A+BK)$ is asymptotically stabilizing. It can be constructed as follows. Let $\Sigma$ denote the positive definite solution of the algebraic Lyapunov equation
\[
(A+BK)^\tr \Sigma (A+BK) + I = \Sigma,
\]
let $|\lambda_\mathrm{max}(A+BK)| \leq \overline \lambda < 1$ be an upper bound on the spectral radius of the matrix $A+BK$, and let $\Lambda > 0$ be the Lipschitz constant of $\Theta$ with respect to the weighted Euclidean norm, $\Vert x \Vert_{\Sigma} \defeq \sqrt{x^\tr \Sigma x}$, such that
\[
\forall x \in X, \qquad \Theta(x,0) \ \leq \ \Lambda \min_{x' \in Y^\star} \| x - x' \|_{\Sigma} \; .
\]
Next, we claim that the function
\[
\forall x \in X, \qquad \Pi(x) \ \defeq \ \frac{\Lambda}{1-\overline \lambda} \left[ \min_{x' \in Y^\star} \| x - x' \|_{\Sigma} \right]
\]
satisfies~\eqref{eq::LMd}. In order to prove this, notice that the inequality
\[
\forall x \in X, \qquad \underset{w \in W}{\max} \; \Pi( f(x,Kx,w) ) \ \leq \ \overline \lambda \cdot \Pi(x) 
\]
holds by construction of $\Sigma$, $\Pi$ and $Y^\star$. Thus, we have
\begin{align}
& \Theta(x,0) + \underset{w \in W}{\max} \; \Pi( f(x,Kx,w) ) \notag \\
&\ \leq \ \Lambda \min_{x' \in Y^\star} \| x - x' \|_{\Sigma} + \overline \lambda \cdot \Pi(x) \notag \\
&\ = \ \left[ \Lambda + \frac{\Lambda \overline \lambda}{1-\overline \lambda} \right] \min_{x' \in Y^\star} \| x - x' \|_{\Sigma} = \Pi(x)
\end{align}
for all $x \in X$; that is, $\Pi$ satisfies~\eqref{eq::LMd}.
Next, the associated ambiguity measure
\begin{eqnarray}
\label{eq::MM}
M(P) = \max_{p \in P} \int \Pi \, \mathrm{d}p
\end{eqnarray}
can be used as an associated terminal cost that satisfies Assumption~\ref{ass::terminal}. As this result is of high practical relevance, we summarize it in the form of the following lemma.

\begin{lemma}
\label{lem::assumptions}
Let $L$ and $M$ be defined as in~\eqref{eq::LL} and~\eqref{eq::MM}. If the functions $\Theta$ and $\Pi$ are non-negative and Lipschitz continuous such that~\eqref{eq::LMd} is satisfies and if $0 \in U_\mathrm{c}$, then $L$ and $M$ satisfy all requirements from Assumptions~\ref{ass::proper} and~\ref{ass::terminal}.
\end{lemma}

\pf
Because $\Theta$ and $\Pi$ are Lipschitz continuous $L$ and $M$ are---by construction---proper ambiguity measures and Assumption~\ref{ass::proper} is satisfied. Moreover, since $\Theta$ and $\Pi$ are non-negative, $L$ and $M$ are non-negative. Next,~\eqref{eq::MM} implies that
\begin{eqnarray}
M(F(P,\mu)) &\overset{\eqref{eq::MM}}{=}& \max_{p \in F(P,\mu)} \int \Pi \, \mathrm{d}p \notag \\[0.16cm]
&=& \max_{p \in P, \omega \in \Omega} \int \hspace{-0.1cm} \int \Pi( f(x,K x,w) ) \, p(\mathrm{d}x) \, \omega( \mathrm{d}w ) \notag \\[0.16cm]
&\leq& \max_{p \in P} \int \max_{w \in W} \Pi( f(x,Kx,w) ) \, p(\mathrm{d}x) \; ,
\label{eq::Auxx}
\end{eqnarray}
where the second equation holds for the offset-free ancillary feedback law $\mu(x) = Kx$. Furthermore, according to~\eqref{eq::LMd}, we have
\begin{eqnarray}
\label{eq::Auxx2}
\max_{w \in W} \; \Pi( f(x,Kx,w) ) &\leq& \Pi(x) - \Theta(x,0)
\end{eqnarray}
for all $x \in X$. Consequently, we can substitute this inequality in~\eqref{eq::Auxx} finding
\begin{eqnarray}
M(F(P,\mu)) &\overset{\eqref{eq::Auxx},\eqref{eq::Auxx2}}{\leq}& \max_{p \in P} \int \left[ \Pi - \Theta(\cdot,0) \right] \, \mathrm{d}p \notag \\[0.16cm]
&\overset{\eqref{eq::MM},\eqref{eq::LL}}{=}& M(x) - L(P,\mu) \notag
\end{eqnarray}
recalling that this holds for the particular feedback law $\mu(x) = Kx$. In other words, because we assume that $0 \in U_\mathrm{c}$, there exists for every $P \in \mathcal K(X)$ a $\mu \in \mathcal U$ for which
\[
L(P,\mu) + M(F(P,\mu)) \; \leq \; M(P)
\]
and the conditions from Assumption~\ref{ass::terminal} are satisfied. This corresponds to the statement of the lemma.
\qed

\begin{remark}
\label{rem::hybrid}
Notice that many articles on stochastic MPC, for example~\cite{Chatterjee2015,Kouvaritakis2016,Mayne2015}, start their construction of the stage cost by assuming that a nominal (non-negative) cost function $l: X \times U \to \mathbb R_+$ is given. For example, in the easiest case, one could consider the least-squares cost
\[
l(x,u) = x^2 + u^2 \; .
\]
In the above context, however, we cannot simply set \mbox{$\Theta = l$}, as Condition~\eqref{eq::LMd} can only be satisfied if we have $\Theta(y^\star,0) = 0$ for all \mbox{$y^\star \in Y^\star$---but} $Y^\star$ is usually not a singleton. However, one can find a Lipschitz continuous function $\Theta$ that approximates the function
\[
\Theta(x,u) \approx \left\{
\begin{array}{ll}
l(x,u) \; \; & \text{if} \; \; x \notin Y^\star \\[0.1cm]
0 & \text{if} \; \; x \in Y^\star
\end{array}
\right.
\]
up to any approximation accuracy such that $\Theta$ coincides with $l$ on the domain $X \setminus Y^\star$ with high precision. In fact, the approximation is in this context only needed for technical reasons, such that $\Theta$ is Lipschitz continuous. Notice that this construction satisfies the requirements of Lemma~\ref{lem::assumptions} and is, as such, fully compatible with our stability analysis framework. Next, we construct the
hybrid feedback law
\[
\widetilde \mu(y) \; \defeq \; \left\{
\begin{array}{ll}
\mu_{\mathrm{MPC}}(y) \; \; & \text{if} \; \; y \notin Y^\star \\[0.1cm]
Ky & \text{if} \; \; y \in Y^\star \; ,
\end{array}
\right.
\]
which simply switches to the ancillary control law \mbox{$x \to Kx$} whenever the current state is already inside the target region $Y^\star$. This construction is compatible with the stability statements from Theorems~\ref{thm::stability} and~\ref{thm::asympstability}, as we modify the closed-loop system only inside the robust control invariant target region. Of course, this is in the understanding that the control gain $K$ is optimized beforehand and that this linear controller leads to a close-to-optimal control performance (with respect to worst-case expected value of the given cost function~$l$) inside the region $Y^\star$---if not, one needs to work with more sophisticated ancillary controllers and redefine $\mathcal U$. The robust MPC controller is in this case only taking care of the case that the current state is in $X \setminus Y^\star$---but in this region $\Theta$ coincides with the given cost function $l$ as desired.
\end{remark}

\subsection{Implementation Details}
Ambiguity tube MPC can be implemented by pre-computing the stage and terminal cost offline. This has the advantage that, in the online phase, a simple convex optimization problem is solved. For this aim, we pre-compute the central sets
\begin{eqnarray}
Z_{k+1} &\defeq& \bigoplus_{i=0}^k (A+BK)^i W_{\mathrm{c}}
\end{eqnarray}
for all $k \in \{ 0,1, \ldots, N-1 \}$ by using standard set computation techniques~\cite{Blanchini2009}. Similarly, by introducing the Markovian kernel
\[
\forall X' \in \mathcal B(X), \quad \Delta[x](X') \overset{\mathrm{def}}{=} \omega( X' \oplus \{ -(A+BK)x \} ) \; ,
\]
we can pre-compute offset-free measures $q_k$ via the Markovian recursion 
\begin{eqnarray}
\forall k \in \mathbb N, \quad q_{k+1} &=& \int \Delta[x] \, q_{k}(\mathrm{d}x) \\[0.16cm]
q_0 &=& \delta_0 \; . \notag
\end{eqnarray}
For example, if $\omega$ denotes a uniform probability measure with compact zonotopic support, the measures $q_k$ can be computed with high precision by using a generalized Lyapunov recursion in combination with a Gram-Charlier expansion~\cite{Villanueva2020}. After this preparation, we can pre-compute Chebyshev representations of the functions
\begin{eqnarray}
S_k(z,u) &=& \max_{z_\mathrm{c} \in Z_k} \int \Theta( z + z_\mathrm{c} + \overline z, u ) \, q_k(\mathrm{d}\overline z) \notag \\[0.16cm]
\text{and} \qquad S_N(z) &=& \max_{z_\mathrm{c} \in Z_N} \int \Pi( z + z_\mathrm{c} + \overline z ) \, q_N(\mathrm{d}\overline z) \notag 
\end{eqnarray}
with high accuracy, as discussed in~\cite{Villanueva2020}, too. Here, the function $\Theta$ and $\Pi$ are constructed as in the previous section. In particular, if $\Theta$ is convex, as in Example~\ref{ex::theta}, $S_k$ is convex. Similarly, if $\Pi$ is convex, $S_N$ is convex. Finally, the associated online optimization problem,
\begin{eqnarray}
\mathcal V(y) \; = &\underset{v,z}{\min} & \sum_{k=0}^{N-1} S_k(z_k,v_k) + S_N(z_N) \notag \\[0.2cm]
\label{eq::MPCimplementation}
&\mathrm{s.t.}& \left\{
\begin{array}{l}
\forall k \in \{ 0, 1, \ldots, N-1 \}, \\[0.16cm]
z_{k+1} = A z_k + B v_k \\[0.16cm]
z_0 = y \; ,
\end{array}
\right.
\end{eqnarray}
can be solved with existing MPC software. The associated ambiguity tube MPC feedback  has then, by construction, the form
\[
\mu_\mathrm{MPC}(y) = v_0^\star(y) \; ,
\]
where $v_0^\star(y)$ denotes the first element of an optimal control input sequence of~\eqref{eq::MPCimplementation} as a function of $y$. This construction can be further refined by implementing the hybrid control law from Remark~\ref{rem::hybrid}.

\begin{figure*}[t]
\centering
\begin{tabular}{|c|c|}
\hline
\begin{minipage}{0.47\textwidth}
\centering
\vspace{0.3cm}
\textbf{Ambiguity Tube MPC: $\; \; \tau = \frac{1}{10}$}\\[0.4cm]
\includegraphics[width=0.84\textwidth]{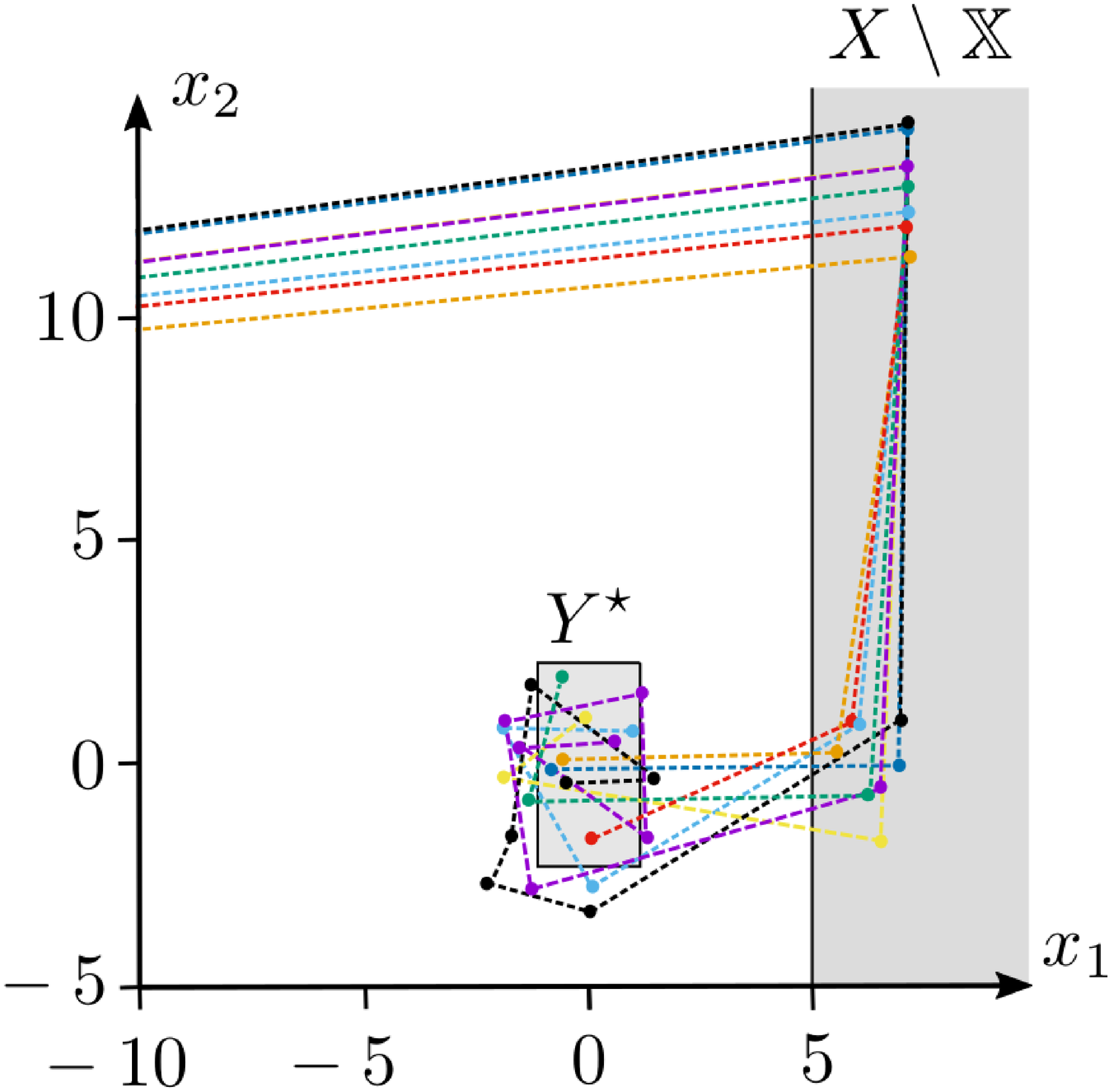}%
\end{minipage}
&
\begin{minipage}{0.47\textwidth}
\centering
\vspace{0.3cm}
\textbf{Ambiguity Tube MPC: $\; \; \tau = 10$}\\[0.4cm]
\includegraphics[width=0.84\textwidth]{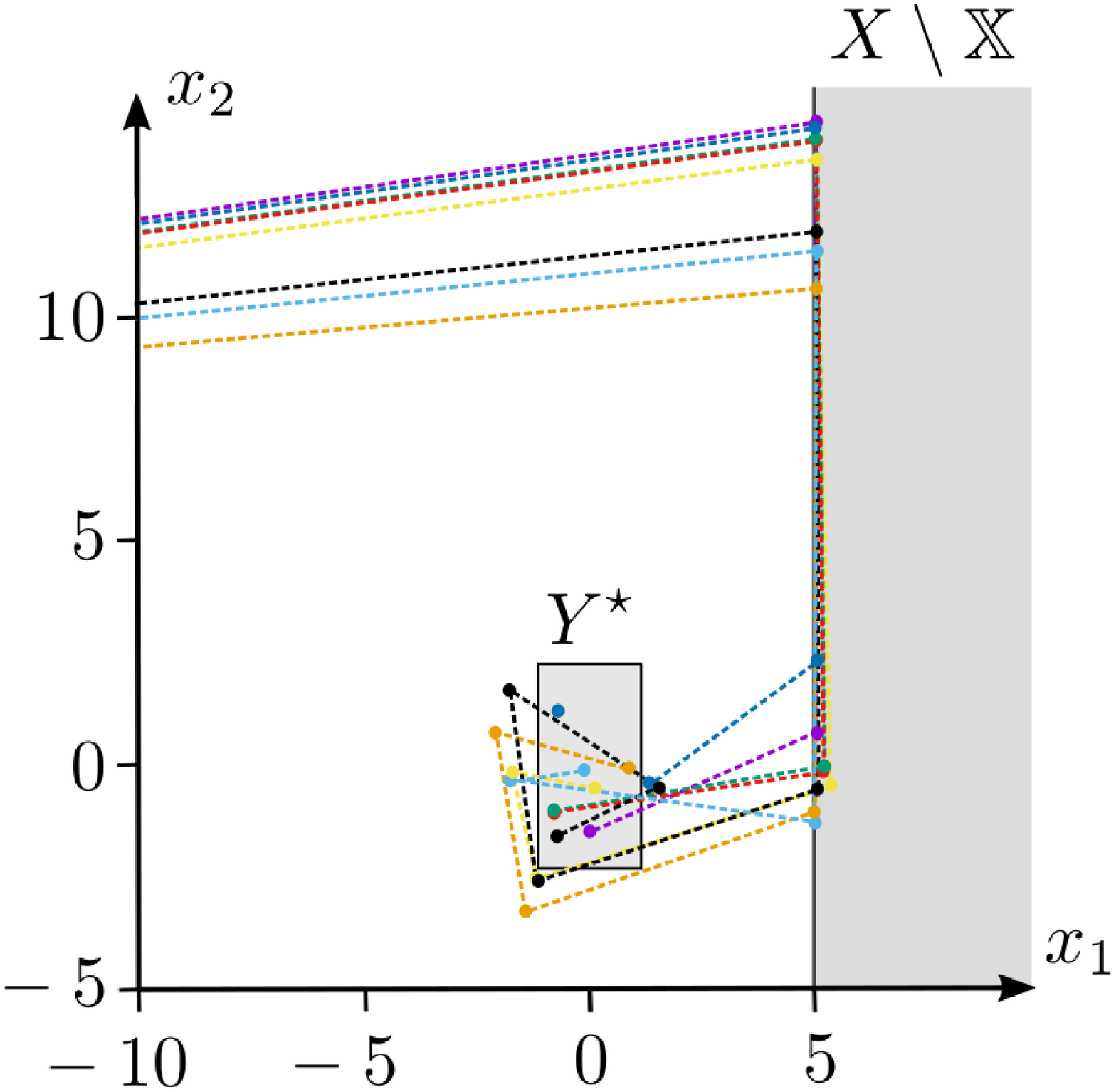}%
\end{minipage}\\
\hline
\end{tabular}
\caption{\label{fig::closedloop} Randomly generated closed-loop scenarios of the ambiguity tube MPC controller~\eqref{eq::MPCimplementation} with (LEFT) $\tau = 10^{-1}$ and (RIGHT) $\tau = 10$. The state constraint, $x_1 \leq 5$, is visualized in the form of the dark gray-shaded infeasible region $x_1 > 5$. For the very small penalty parameter $\tau = \frac{1}{10}$ the controller risks marginal constraint violation for the sake of better nominal control performance. All trajectories converge to the light gray-shaded target region $Y^\star$. Notice that the initial state, $y_0 = [-60,5]^\tr$, would be far off to the left of the figure and it is therefore not visualized---the colored dots correspond to the discrete-time states after the first MPC iteration.} 
 \end{figure*}

\subsection{Numerical Illustration}
This section illustrates the performance of the above ambiguity tube MPC controller for the nonlinear system
\begin{eqnarray}
x_1^+ &=& -\frac{x_1}{8} + x_2 + u_1 + w_1, \notag \\
x_2^+ &=& -\frac{x_1}{2} + \frac{x_2}{4} + u_2 + \cos(x_1) \sin(5 x_2)^3 + w_2 \; . \notag
\end{eqnarray}
In order to write this system in the above form, we introduce the notation
\[
A \, \defeq \, \left(
\begin{array}{rr}
-\frac{1}{8} & 1 \\[0.16cm]
-\frac{1}{2} & \frac{1}{4}
\end{array}
\right) \; , \quad B \, \defeq \, I \quad \text{and} \quad K \, \defeq \, \left(
\begin{array}{rr}
\frac{1}{8} & -\frac{1}{2} \\[0.16cm]
\frac{1}{4} & -\frac{1}{4}
\end{array}
\right)
\]
to denote the system matrices and a suitable ancillary control gain. Notice that the influence of the nonlinear term, $\cos(x_1) \sin(5x_2)^3$, can be over-estimated by introducing the central set $W_\mathrm{c} \defeq \{ 0 \} \times [-1,1]$ (see Remark~\ref{rem::nonlinear}). Additionally, we set $U_\mathrm{c} \defeq [-10,10]^2$. Moreover, the uncertain input is modeled by the distribution measure
\[
\omega(A) \defeq \int_A \rho(w)  \mathrm{d}w
\]
with Radon-Nikodyn derivative (density function)
\[
\rho(w) \; \defeq \; \left\{
\begin{array}{ll}
\mathfrak{d}(w_1) & \text{if} \; w_2 \in [-1, 1] \\
0 & \text{otherwise}
\end{array}
\right\} \; ,
\]
where $\mathfrak{d} \defeq \partial \delta / \partial w$ denotes the standard Dirac distribution. The objective is constructed as in Example~\ref{ex::theta},
\[
\Theta(x,u) \; \defeq \; \| u \|_2^2 + \mathrm{dist}_{2}(x,Y^\star)^2 + \tau \cdot \mathrm{dist}_{1}(x,\mathbb X) \; ,
\]
where $\tau > 0$ is a risk-parameter that is associated with the given state constraint set
\[
\mathbb X \; \defeq \; \left\{ \; x \in \mathbb R^2  \; \middle| \; x_1 \leq 5 \; \right\} \; .
\]
Here, it is not difficult to check that the function
\[
\Pi(x) \; \defeq \; \left[ 2 + \frac{7}{27} \cdot \tau \right] \mathrm{dist}_2(x,Y^\star)^2
\]
satisfies the requirements from Lemma~\ref{lem::assumptions}. Consequently, the associated ambiguity measures $L$ and $M$ satisfy all technical requirements of Theorem~\ref{thm::asympstability}. Notice that the functions $S_k$ in~\eqref{eq::MPCimplementation} are in our implementation pre-computed with high precision such that the convex optimization problem~\eqref{eq::MPCimplementation} can be solved in much less than $1 \, \mathrm{ms}$ by using ACADO Toolkit~\cite{Houska2011}. The prediction horizon of the MPC controller is set to $N = 5$.

Figure~\ref{fig::closedloop} shows two ambiguity tube MPC closed-loop simulations that are both started at the initial point $y_0 = [-60, 5]^\tr$. In the left figure, we have set $\tau = \frac{1}{10}$, which means that the constraint violation penalty is small compared to the nominal control performance objective. Consequently, during the randomly generated closed-loop scenarios marginal constraint violation can be observed. This is in contrast to the right part of Figure~\ref{fig::closedloop}, which shows randomly generated closed-loop scenarios for the case $\tau = 10$, leading to much smaller expected constraint violations at risk. In all cases, that is, independent of how the penalty parameter $\tau > 0$ is chosen and independent of the particular uncertainty scenario, the closed-loop trajectories converge to the terminal region $Y^\star$ after a short transition period. In this particular example, we observe that this happens typically after $3$ to $8$ discrete-time steps, which confirms the robust asymptotic convergence statement of Theorem~\ref{thm::asympstability}.

\section{Conclusion}
\label{sec::conclusions}
This paper has presented a coherent measure-theoretic framework for analyzing the stability of a rather general class of ambiguity tube MPC controllers. In detail, we have proposed a Wasserstein-Hausdorff metric leading to our first main result in Theorem~\ref{thm::existence}, where conditions for the existence of a continuous value function of ambiguity tube MPC controllers have been established. Moreover, Theorem~\ref{thm::martingale} has built upon this topological framework to establish conditions under which the stage and terminal cost are proper ambiguity measures, such that the cost function of the MPC controller can be turned into a non-negative supermartingale along the trajectories of the stochastic closed-loop system. Related stochastic stability and convergence results for ambiguity tube MPC have been summarized in Theorems~\ref{thm::stability} and~\ref{thm::asympstability}. 

In the sense that Lemma~\ref{lem::assumptions} proposes a practical strategy for constructing stabilizing terminal costs for stochastic and ambiguity tube MPC, the current article has outlined a path towards a more consistent stability theory that goes much beyond the existing convergence results from~\cite{Chatterjee2015,Kouvaritakis2016,Munoz2020}. At the same time, however, it should also be pointed out that these results are based on a slightly different strategy of modeling the stage cost of the MPC controller, as discussed in Remark~\ref{rem::hybrid}, where it is also explained why it may be advisable to use a hybrid feedback control law that switches to a pre-optimized ancillary controller whenever the state is inside its associated target region $Y^\star$.

Last but not least, as much as this article has been attempting to make a step forward, towards a more consistent stability theory and practical formulation of robust MPC, it should also be stated that many open problems and conceptual challenges remain. In the line of this paper, for instance, a discussion of more advanced representations of ambiguity set representations and handling of nonlinearities, a more in-depth analysis of the interplay of the choice of $\mathcal U$, the performance of the ambiguity tube controller, and its computational tractability, bounds on the concentrations of the state distributions, economic objectives for ambiguity tube MPC, as well as a deeper analysis of risk measures and related issues of recursive feasibility are only a small and incomplete selection of open problems in the field of distributionally robust MPC.

\bibliographystyle{plain}
\bibliography{references}

\end{document}